\documentclass[12pt, a4paper]{article}
\usepackage{amssymb,latexsym, amsmath,amscd,amsfonts,amsthm,mathbbol}
\newtheorem{prop}{Proposition}[section]
\newtheorem{defi}[prop]{Definition}

\newtheorem{lem}[prop]{Lemma}
\newtheorem{thm}[prop]{Theorem}

\newtheorem{remar}[prop]{Remark}
\newtheorem{cor}[prop]{Corollary}
\newtheorem{nota}[prop]{Notation}
\DeclareMathOperator{\every}{\forall}
\DeclareMathOperator{\Aut}{Aut}
\DeclareMathOperator{\Hom}{Hom}
\DeclareMathOperator{\Ind}{Ind}
\DeclareMathOperator{\cInd}{c-Ind}
\DeclareMathOperator{\Res}{\mid}

\DeclareMathOperator{\GL}{GL}
\DeclareMathOperator{\Ker}{Ker}
\DeclareMathOperator{\WD}{WD}
\DeclareMathOperator{\Gal}{Gal}

\newcommand{\Rest}[3]{#3\Res_{#2} }
\newcommand{\Indu}[3]{\Ind_{#1}^{#2}#3}
\newcommand{\cIndu}[3]{\cInd_{#1}^{#2}#3}
\newcommand{\cgate}[2]{#2^{#1}}
\newcommand{\tcap}[3]{#2 \cap \cgate{#1}{#3}} 
\newcommand{\IndRes}[4]{\Indu{#1}{#2}{\Rest{#3}{#1}{#4}}}
\DeclareMathOperator{\cRes}{Res}

\newcommand{\cRest}[2]{\cRes^G_#1 #2}
\newcommand{\IndResg}[5]{\IndRes{{ \tcap{#1}{#2}{#4}}}{#2}{{\cgate{#1}{#4}}}{#5^#1}} 
\newcommand{\Mplus}[4]{\bigoplus_{{#1 \in #2\backslash#3/#4}}}
\newcommand{\mackey}[5]{\Mplus{#1}{#2}{#3}{#4} \IndResg{#1}{#2}{#3}{#4}{#5}}

\DeclareMathOperator{\End}{End}

\newcommand{\CC}{{\mathbb C}}
\newcommand{\RR}{{\mathbb R}}
\newcommand{\QQ}{{\mathbb Q}}
\newcommand{\ZZ}{{\mathbb Z}}
\newcommand{\NN}{{\mathbb N}}

\newcommand{\pp}{{\mathfrak p}}
\newcommand{\FF}{{\mathbb F}}

\newcommand{\KK}{{\mathcal K}}
\newcommand{\LL}[1]{\mathcal L(#1)}

\newcommand{\MM}{{\mathbf {M}}}
\newcommand{\MMe}[1]{\MM(\frac{N}{e},#1)}

\newcommand{\curlyu}{{\mathfrak A}}
\newcommand{\uu}[1]{\mathbf U^{#1}(\curlyu)}
\newcommand{\uz}{\mathbf U(\curlyu)}
\newcommand{\urlyb}{\mathfrak B}
\newcommand{\curlyb}[1]{\mathfrak B_{#1}}

\newcommand{\ub}[2]{\mathbf U^{#1}(\curlyb{#2})}
\newcommand{\ubz}[1]{\mathbf U(\curlyb{#1})}
\newcommand{\Ip}{\mathfrak I(\pi)}

\newcommand{\curlyp}{\mathfrak P}
\newcommand{\curlypp}[1]{\curlyp^{#1}}
\newcommand{\curlyq}[1]{\mathfrak Q_{#1}}
\newcommand{\curlyqq}[2]{\curlyq{#2}^{#1}}
\newcommand{\oF}{\mathfrak o_F}
\newcommand{\pF}{\mathfrak p_F}
\newcommand{\oE}{\mathfrak o_E}
\newcommand{\pE}{\mathfrak p_E}
\newcommand{\kF}{\mathfrak k_F}
\newcommand{\kE}{\mathfrak k_E}
\newcommand{\kU}{\mathfrak K(\curlyu)}
\newcommand{\EJ}{E^{\times}J}
\newcommand{\FK}{F^{\times}K}
\newcommand{\GLNe}{\GL_{\frac{N}{e}}(\oF)}
\newcommand{\IG}[3]{I_G(#1,#2| #3)}

\newcommand{\IB}[4]{I_{B^{\times}}(#2,#3| #4)}
\newcommand{\brac}[3]{\langle #1,#2 \rangle_{#3}}
\newcommand{\uni}[1]{\pi_F^{{\alpha_#1}}}
\newcommand{\al}[1]{\alpha_#1}
\newcommand{\eUof}{e(\curlyu|\oF)}
\newcommand{\eBoe}[1]{e(\curlyb{#1}|\oE)}

\newcommand{\JbU}{J(\beta,\curlyu)}
\newcommand{\JibU}[1]{J^{#1}(\beta,\curlyu)}
\newcommand{\HibU}[1]{H^{#1}(\beta,\curlyu)}
\newcommand{\HigU}[1]{H^{#1}(\gamma)}
\newcommand{\HitU}[1]{H^{#1}(\beta)}
\DeclareMathOperator{\tr}{tr}
\newcommand{\Kg}{\mathcal K(g)}
\newcommand{\XX}{\mathcal X}

\newcommand{\psa}{\psi_A}
\newcommand{\psb}[1]{\psi_{#1}}

\newcommand{\valu}{\nu_{\curlyu}}
\newcommand{\oo}{\mathfrak o}
\newcommand{\CUm}[2]{\mathcal C(\curlyu,#1,#2)}
\newcommand{\LLL}{\mathcal {L}}
\title{Unicity of types for supercuspidals}
\author{Vytautas Paskunas}
\date{\today.}
\begin{document} 
\maketitle
\tableofcontents
\pagebreak
\section{Introduction}
Let $F$ be a non-Archimedean local field with a finite residue field $\kF$. Let $\oF$ be its complete discrete valuation ring, $\pF$ the maximal ideal of
 $\oF$, and $q_F$ the size of $\kF$. Moreover, let $N > 1$, $G=\GL_N(F)$ and $K=\GL_N(\oF)$. Further, let $W_F$ be the Weil group of $F$ and $I_F$ be the
 inertia group of $F$. All the representations considered in this paper are over $\CC$.

\begin{defi} Let $\pi$ be a smooth irreducible supercuspidal representation of $G$, then we define the \textbf{inertial support} $\Ip$ of $\pi$ to be:
$$\Ip=\{\pi':\pi'\cong  \pi \otimes \chi \circ \det\}$$
where $\chi$ is some unramified quasicharacter of $F^{\times}$.
\end{defi} 
\begin{defi} Suppose $H$ is a compact open subgroup of $G$, $\tau$ a smooth irreducible representation of $H$ and $\pi$ a smooth irreducible supercuspidal representation of $G$, then $(H,\tau)$ is a \textbf{type} for $\Ip$, if for all smooth irreducible representations $\pi'$ of $G$:
$$\Rest{G}{H}{\pi'}\mbox{ contains }\tau \Leftrightarrow \pi'\in\Ip$$
where $\chi$ is some unramified quasicharacter of $F^{\times}$.
\end{defi}
Our main result is:
\begin{thm}\label{main} Let $\pi$ be a smooth irreducible supercuspidal representation of $G$, then there exists a smooth irreducible representation $\tau$ of $K$ depending on $\Ip$, such that $(K,\tau)$ is a type for $\Ip$. Moreover, $\tau$ is  unique (up to isomorphism) and it occurs in $\Rest{G}{K}{\pi}$ with multiplicity one. 
\end{thm}
This implies a kind of inertial local Langlands correspondence:
\begin{cor}\label{supercusp}Let $\varphi$ be a smooth $N$-dimensional representation of $I_F$, which extends to a smooth \underline{irreducible} Frobenius semisimple representation of $W_F$, then there exists a unique (up to isomorphism) smooth irreducible representation $\tau(\varphi)$ of $K$, such that for any smooth irreducible infinite dimensional representation $\pi$ of $G$:
$$\Rest{G}{K}{\pi}\mbox{ contains }\tau(\varphi)\Leftrightarrow \Rest{{W_F}}{{I_F}}{\WD(\pi)}\cong\varphi$$
where $\WD(\pi)$ is a Weil-Deligne representation of $W_F$ corresponding to $\pi$ via the local Langlands correspondence.
\end{cor} 

Our result and methods generalise the case, when $N=2$, which was considered by G.\,Henniart in \cite{ap}. The paper heavily relies on the classification of supercuspidals due to C.\,Bushnell and P.\,Kutzko. The existence of such $\tau$ is almost immediate from \cite{bk}(6.2.3), the difficult part is proving uniqueness. 

The paper is structured as follows. We recall some facts and definitions from Bushnell-Kutzko theory in sections \ref{her}-\ref{simt}. In section \ref{setup} we introduce some of our own notation. From Bushnell-Kutzko theory we know that every supercuspidal representation is induced from an open compact-mod-centre subgroup of $G$. The restriction of $\pi$ to $K$ results in a Mackey's decomposition. The representation coming from the double coset, which contains $1$, is our $\tau$. We prove that in section \ref{existence}. Then in section \ref{important} we prove that under certain conditions an irreducible summand of $\Rest{G}{K}{\pi}$ cannot be a type. In section \ref{represent} we choose a nice representative from each double coset. Then we have to consider two different cases, namely sections \ref{a} and \ref{b}. The idea is that unless the double coset contains $1$, then any irreducible summand of the representation coming from the double coset in the Mackey's decomposition of $\Rest{G}{K}{\pi}$ cannot be a type by section \ref{important}. Finally in section \ref{over} we prove the main result. 

This is a part of my PhD thesis at the University of Nottingham.
 I would like to thank my supervisor M.\,Spie\ss, and also C.\,Bushnell and G.\,Henniart for  
their comments on the earlier draft.
 
\section{Notation and Preliminaries}
Let $F$ be a non-Archimedean local field with a finite residue field $\kF$. Let $\oF$ be its complete discrete valuation ring, $\pF$ the maximal ideal of
 $\oF$, and $q_F$ the size of $\kF$. 
 Moreover, let $V$ be an $F$-vector space of finite dimension $N> 1$, $A=\End_F(V)$, $G=\Aut_F(V)$. Further, let $\psi_F$ be a fixed continuous additive character of the group $F$, with conductor $\pF$, and let $$\psi_A=\psi_F\circ\tr_{A/F}.$$ 
\subsection{Intertwining}
If $H$ is a subgroup of $G$ and $\rho$, $\tau$ are representations of $H$, let $$\brac{\rho}{\tau}{H}=\dim_{\CC}\Hom_H(\rho,\tau).$$ If $g\in G$, then let $$H^g=gHg^{-1}$$ and $\rho^g$, be a representation of $H^g$, $$\rho^g(x)=\rho(g^{-1}xg)\mbox{, } \forall x\in H^g.$$
We say  $g$ \emph{intertwines} $\tau$ and $\rho$ in $G$, if $$\brac{\tau}{\rho^g}{H\cap H^g}\neq 0.$$ The set of all $g\in G$, which intertwine $\tau$ and $\rho$ is called the \emph{intertwining} of $\tau$ and $\rho$ and is denoted by $\IG{\tau}{\rho}{H}$. 
 
\subsection{Hereditary orders}\label{her}
For a complete account of hereditary orders we refer the reader to  \cite{bk}\S1, \cite{bu1} and \cite{bf}. Everything below is taken from \cite{bk}\S1.1.

Let $\curlyu$ be an $\oF$-order in $A$, then $\curlyu$ is (left) \emph{hereditary} if every (left) $\curlyu$-lattice is $\curlyu$-projective. 

\subsubsection{Lattice chains}
An $\oF$-\emph{lattice chain} $\mathcal L$ in $V$ is a sequence $\{L_i:i\in \ZZ\}$, such that 
\begin{itemize}
\item[(i)]$L_{i+1}\subsetneq L_i$, $i\in\ZZ$ 
\item[(ii)] there exists $e\in \ZZ$ such that $\pi_F L_i=L_{i+e}$ for all $i\in\ZZ$. 
\end{itemize}
The integer $e=e(\mathcal L)$ is uniquely determined and is called an $\oF$-\emph{period} of  $\mathcal L$.
\subsubsection{Hereditary orders}
Hereditary orders in $A$ are in bijection with lattice chains in $V$. 
Given a lattice chain $\mathcal L$  we define:
$$\End^n_{\oF}(\LLL)=\{ x\in A:xL_i\subseteq L_{i+n}, i\in\ZZ \}$$
for each $n\in\ZZ$. Then $\curlyu=\curlyu(\LLL)=\End^0_{\oF}(\LLL)$ is a hereditary $\oF$-order in $A$. We can recover $\LLL$ from $\curlyu$ up to a shift in the index: $\LLL$ is the set of all $\oF$-lattices in $V$, which are $\curlyu$ modules. The lattices $\End^n_{\oF}(\LLL)$ are $(\curlyu,\curlyu)$-bimodules. Let $\curlyp$ be the Jacobson radical of $\curlyu$, then $$\curlyp=\End^1_{\oF}(\LLL).$$ As a fractional ideal of $\curlyu$, the radical $\curlyp$ is invertible, and we have:
$$\curlyp^n=\End^n_{\oF}(\LLL)\mbox{, } n \in \ZZ.$$
In particular,
$$\curlyp^n L_i=L_{i+n}\mbox{, } i, n\in \ZZ.$$
The $\oF$-period $e$ of $\LLL$ is also a function of $\curlyu=\curlyu(\LLL)$, given by
$$\pF\curlyu=\curlyp^{e}.$$ 
So we will write $e=e(\LLL)=\eUof$. We define a sequence of compact open subgroups in $G$:
$$\uu{0}=\uz=\curlyu^{\times},$$
$$\uu{n}=1+\curlypp{n}\mbox{, }n\ge 1.$$
Also, for $\curlyu=\curlyu(\LLL)$, we set
$$\kU=\{x\in G:xL_i\in \mathcal L,i\in\ZZ\}=\{x\in G:x^{-1}\curlyu x=\curlyu\}.$$
This is an open compact-mod-centre subgroup of $G$, and the $\uu{n}$, for $n\ge 0$, are normal subgroups of it. In particular, $\uz$ is the unique maximal compact open subgroup of $\kU$.  

There is also a "valuation" map associated with the hereditary order $\curlyu$. Define:
$$\valu(x)=\max\{n\in\ZZ: x\in \curlyp^n\}\mbox{, }x\in A,$$
with the understanding that $\valu(0)=\infty$. In particular, if $x\in\kU$, then $\valu(x)=n$, where
$$x\curlyu=\curlyu x=\curlyp^n.$$
This induces an exact sequence:
\[\begin{CD}{1}@>>>{\uz}@>>>{\kU}@>{\valu}>>{\ZZ}  \end{CD} \]
Given $\curlyu=\curlyu(\LLL)$ we have a canonical isomorphism
$$\curlyu/\curlyp\cong\prod^{e-1}_{i=0} \End_{\kF}(L_i/L_{i+1}),$$ where $e=\eUof$. And, we can always choose a basis for $V$, such that $\curlyu$ is identified with block upper triangular matrices modulo $\pF$. We say $\curlyu$ is \textit{principal} if $(L_i:L_{i+1})=(L_0:L_1)$, for all $i$. In that case, every block has size $\frac{N}{e}\times\frac{N}{e}$.
 
\subsubsection{The character $\psb{b}$}

Let $n$ and $m$ be integers satisfying 
$n>m\ge [\frac{n}{2}]\ge 0,$
where $[x]$ denotes the greatest integer $\le x$, for $x\in \RR$. We then have a canonical isomorphism 
\[\begin{CD}{\uu{m+1}/\uu{n+1}}@>{\cong}>>{\curlyp^{m+1}/\curlyp^{n+1}},  \end{CD} \]
given by $x\mapsto x-1$. This leads to an isomorphism 
\[\begin{CD}{(\uu{m+1}/\uu{n+1})^{\wedge}}@>{\cong}>>{\curlyp^{-n}/\curlyp^{-m}},  \end{CD} \] where "hat" $^{\wedge}$ denotes Pontryagin dual. Explicitly, this is given by 
$$b+\curlyp^{-m}\mapsto \psb{A,b}=\psb{b}\mbox{, } b\in \curlyp^{-n}\mbox{, where}$$
$$\psb{b}(1+x)=\psb{A}(bx)\mbox{, }x\in\curlyp^{m+1}.$$
 
\subsection{Strata}

For details we refer the reader to \cite{bk}\S1.

A \emph{stratum} is a 4-tuple $[\curlyu,n,m,b]$ consisting of a hereditary order $\curlyu$, integers $n>m$, and $b\in A$, such that $\valu(b)\ge -n$. 

\subsubsection{Equivalence}
We define an equivalence relation on the set of strata: $$[\curlyu_1,n_1,m_1,b_1]\sim [\curlyu_2,n_2,m_2,b_2]\mbox{, if}$$ $$b_1+\curlyp_1^{-m_1}=b_2+\curlyp_2^{-m_2}.$$ Equivalence implies  $\curlyu_1=\curlyu_2=\curlyu$, $m_1=m_2$. Moreover,  if $n_1=-\valu(b_1)$ and $n_2=-\valu(b_2)$, then $n_1=n_2$,
 see \cite{bk}(1.5.2).

\subsubsection{Simple strata}
A stratum $[\curlyu,n,m,\beta]$ is \emph{pure} if 
\begin{itemize}
\item[(i)]the algebra $E=F[\beta]$ is a field,
\item[(ii)] $E^{\times}\subset\kU$,
\item[(iii)] $\valu(\beta)=-n$. 
\end{itemize}
It is called \emph{simple} if, in addition 
\begin{itemize}
\item[(iv)] $m<-k_0(\beta,\curlyu)$.
\end{itemize} 
The definition of $k_0(\beta,\curlyu)$ is rather technical, so we refer the reader to \cite{bk}(1.4.5). We will need to know that $k_0(\beta,\curlyu)$ is an integer and $$k_0(\beta,\curlyu)\ge\valu(\beta).$$ 

Suppose $[\curlyu,n,m,\beta]$ is a simple stratum, then we define 
$$B_{\beta}=\{x\in A:\beta x=x\beta\}=\End_E(V).$$
Let $$\curlyb{\beta}=\curlyu\cap B_{\beta}.$$
Since $E^{\times}\subset\kU$, we can view $\mathcal L$ as an $\oE$ lattice chain. Hence  $\curlyb{\beta}$ is a hereditary order in $B_{\beta}$. We define $$\curlyqq{n}{\beta}=\curlypp{n}\cap B_{\beta},$$
 $$\ubz{\beta}=\uz\cap B_{\beta}$$ and $$\ub{n}{\beta}=\uu{n}\cap B_{\beta}.$$ All the notions above coincide with the ones defined for $\curlyu$. We also have $$\eBoe{\beta} e(E|F)=\eUof,$$ where $e(E|F)$ is a ramification index of $E$ over $F$, since $\pi_E^{e(E|F)}.L_i=\pi_F.L_i$, for all $L_i\in \mathcal{L}$. 
\subsubsection{Tame corestriction}
Let $E/F$ be a subfield of $A$, with the centraliser $B$. 
A \emph{tame corestriction} on $A$ relative to $E/F$, is a $(B,B)$-bimodule homomorphism, $s:A\rightarrow B$, such that $s(\curlyu)=\curlyu\cap B$ for every hereditary $\oF$ order $\curlyu$ in $A$, which is normalised by $E^{\times}$. We will need the following result \cite{bk}(1.3.4).

Let $\psi_E$ be a continuous additive character of $E$, with conductor $\pE$ and let $\psi_B=\psi_E\circ\tr_{B/E}$, then there exists a unique map $s:A\rightarrow B$, such that $$\psi_A(ab)=\psi_B(s(a)b)\mbox{, }\forall a\in A\mbox{, } \forall b\in B.$$ This map is a tame corestriction relative to $E/F$.
\subsubsection{Approximation of simple strata}
We will use the following result \cite{bk}(2.4.1).
\begin{enumerate}
\item[(i)] Let $[\curlyu,n,m,\beta]$ be a pure stratum in $A$, then there exists a simple stratum $[\curlyu,n,m,\gamma]$ in $A$, such that
$$ [\curlyu,n,m,\beta]\sim[\curlyu,n,m,\gamma].$$
Among all the pure stratum $[\curlyu,n,m,\beta']$  equivalent to $[\curlyu,n,m,\beta]$ the simple ones are precisely those for which the field extension $F[\beta']/F$ has minimal degree.
\item[(ii)]Let $[\curlyu,n,m,\gamma_1]$, $[\curlyu,n,m,\gamma_2]$ be simple strata in $A$, which are equivalent to each other, then $$k_0(\gamma_1,\curlyu)=k_0(\gamma_2,\curlyu).$$  
\item[(iii)] Let $[\curlyu,n,r,\beta]$ be a pure stratum in $A$, with $r=-k_0(\beta,\curlyu)$. Let $[\curlyu,n,r,\gamma]$ be a simple stratum in $A$ which is equivalent to $[\curlyu,n,r,\beta]$, let $s_{\gamma}$ be a tame corestriction on $A$ relative to $F[\gamma]/F$, let $B_{\gamma}$ be the $A$-centraliser of $\gamma$, and $\curlyb{\gamma}=\curlyu\cap B_{\gamma}$. Then $[\curlyb{\gamma},r,r-1,s_{\gamma}(\beta-\gamma)]$ is equivalent to a simple stratum in $B_{\gamma}$.
\end{enumerate}

We will also need \cite{bk}(2.2.8).

Let $[\curlyu,n,m,\beta]$ be a simple stratum in $A$. Let $B$ be the $A$-centraliser of $E=F[\beta]$, and $\urlyb=B\cap\curlyu$. Let $b\in A$ with $\valu(\beta)=-r$, and let $s$ be a tame corestriction on $A$ relative to $F[\beta]/F$. Suppose that the stratum $[\urlyb,m,m-1,s(b)]$ is equivalent to some simple stratum $[\urlyb,m,m-1,c]$ in $B$. Then $[\curlyu,n,m-1,\beta+b]$ is equivalent to a simple stratum $[\curlyu,n,m-1,\beta_1]$. Moreover, if $E_1=F[\beta, c]$, $K=F[\beta_1]$, we have 
\begin{itemize}
\item[(i)]$e(K|F)=e(E_1|F)$, $f(K|F)=f(E_1|F)$;
\item[(ii)]$k_0(\beta_1,\curlyu)= \max\{k_0(\beta,\curlyu)$, $k_0(c,\urlyb)\}$.
\end{itemize}

\subsection{Simple types}\label{simt}

Let $[\curlyu,n,0,\beta]$ be a simple stratum, and let $r=-k_0(\beta,\curlyu)$. 

\subsubsection{Groups $J(\beta,\curlyu)$ and $H(\beta, \curlyu)$}
To a simple stratum $[\curlyu,n,0,\beta]$ we can associate compact open subgroups of $\uz$: $J(\beta,\curlyu)$ and $H(\beta, \curlyu)$, see \cite{bk}(3.1.14). Both of them have natural filtrations by normal subgroups: $$J^m(\beta,\curlyu)=J(\beta,\curlyu)\cap\uu{m},$$ $$H^m(\beta, \curlyu)=H(\beta, \curlyu)\cap\uu{m}.$$
The groups $J^m(\beta,\curlyu)$ and $H^m(\beta, \curlyu)$ are normalised by $\mathfrak {K}(\curlyb{\beta})$, for all $m\ge 0$.
Moreover, $H(\beta, \curlyu)$ is a subgroup of $J(\beta,\curlyu)$ and $H^m(\beta, \curlyu)$ are normal in $J(\beta,\curlyu)$, for $m\ge 1$. We will drop various indices, when the meaning is clear. 

We have the following decompositions: for $0\le m\le [\frac{r}{2}]+1$,
$$H^m(\beta, \curlyu)=\ub{m}{\beta}H^{[\frac{r}{2}]+1}(\beta, \curlyu)$$  
and for $0\le m\le [\frac{r+1}{2}]$
$$J^m(\beta, \curlyu)=\ub{m}{\beta}J^{[\frac{r+1}{2}]}(\beta, \curlyu)$$         where square brackets denote the integer part, see \cite{bk}(3.1.15). 

\subsubsection{Simple characters $\CUm{m}{\beta}$}  
We can define a very special set of linear characters $\CUm{m}{\beta}$ of $H^{m+1}(\beta, \curlyu)$, called \emph{simple characters}, see \cite{bk}(3.2). We will need the following properties:

For $0\le m\le [\frac{r}{2}]$ the restriction induces a surjective map $$\CUm{m}{\beta}\rightarrow \CUm{[\frac{r}{2}]}{\beta}$$
The fibres of this map are of the form $\theta.X$, where $\theta\in\CUm{m}{\beta}$ and $X$ is the group of characters of $\ub{m+1}{\beta}/\ub{[\frac{r}{2}]}{\beta}$, which factor through the determinant $\det_{B_{\beta}}$, see \cite{bk}(3.2.5).

If $n=1$, then $H^1(\beta,\curlyu)=J^1(\beta,\curlyu)=\uu{1}$ and $\CUm{0}{\beta}=\{\psb{\beta}\}$, see \cite{bk}(3.1.7) and (3.2.1).

\subsubsection{Intertwining of simple characters} For $0\le m\le [\frac{r}{2}]$ and $\theta\in\CUm{m}{\beta}$, we have 
$$\IG{\theta}{\theta}{H^{m+1}(\beta,\curlyu)}=J^1(\beta,\curlyu)B_{\beta}^{\times}J^1(\beta,\curlyu)$$
see \cite{bk}(3.3.2).

For $i=1,2$, let $[\curlyu,n,m,\beta_i]$ be simple strata with $m\ge 0$. Suppose there exists $\theta_i\in \CUm{m}{\beta_i}$, which intertwine in $G$. Then there exists $x\in\uz$ such that $$\CUm{m}{\beta_1}=\CUm{m}{x^{-1}\beta_2x}$$
and conjugation by $x$ carries $\theta_1$ to $\theta_2$, see \cite{bk}(3.5.11).
 
If $n>1$ and $r=1$ and $\theta\in \CUm{0}{\beta}$, then there exists a simple stratum $[\curlyu,n,1,\gamma]$, such that $[\curlyu,n,1,\beta]\sim[\curlyu,n,1,\gamma]$, $H^1(\beta,\curlyu)=H^1(\gamma,\curlyu)$ and  $$\theta=\theta_0\psb{c}$$ where $\theta_0\in\CUm{0}{\gamma}$ and $c=\beta-\gamma$, see \cite{bk}(3.2.3).

Moreover, $\IG{\theta}{\theta_0}{H^1(\gamma,\curlyu)}=\emptyset$, see \cite{bk}(3.5.12).
 
\subsubsection{Representations $\eta$ and $\kappa$}
If $\theta\in\CUm{0}{\beta}$, then there exists a unique irreducible smooth representation $\eta$ of $J^1(\beta, \curlyu)$, such that $\Rest{J^1}{H^1}{\eta}$ contains $\theta$, see \cite{bk}(5.1.1).

Given $\eta$, there exists a smooth irreducible representation $\kappa$ of $J(\beta, \curlyu)$ , such that $\Rest{J}{J^1}{\kappa}\cong\eta$ and $B_{\beta}^{\times}\subset \IG{\kappa}{\kappa}{J}$. We say $\kappa$ is a $\beta$-\emph{extension} of $\eta$, see \cite{bk}(5.2.2).

\subsubsection{Simple types}
A \emph{simple type} in $G$ is one of the following \cite{bk}(5.5.10):

1. An irreducible representation $\lambda=\kappa\otimes\sigma$ of $J=J(\beta, \curlyu)$, where $\curlyu$ is a principal $\oF$ order in $A$,  $[\curlyu,n,0,\beta]$ is a simple stratum. For some $\theta\in\CUm{0}{\beta}$, $\kappa$ is a $\beta$-extension of $\eta$, the unique irreducible representation of $J^1(\beta,\curlyu)$ containing $\theta$. Let $E=F[\beta]$, then 
$$J(\beta, \curlyu)/J^1(\beta, \curlyu)\cong\ubz{\beta}/\ub{1}{\beta}\cong\GL_f(\kE)^e$$ for some integers $e$ and $f$ and $\sigma$ is a lift of representation $\sigma_0\otimes\ldots\otimes\sigma_0$, where $\sigma_0$ is an irreducible cuspidal representation of $GL_f(\kE)$.

2. An irreducible representation $\sigma$ of $\uz$, where $\curlyu$ is a principal $\oF$ order in $A$. We have $\uz/\uu{1}\cong\GL_f(\kF)^e$, for some integers $e$ and $f$. Then $\sigma$ is a lift of $\sigma_0\otimes\ldots\otimes\sigma_0$, where $\sigma_0$ is an irreducible cuspidal representation of $GL_f(\kF)$.

The second part, can be viewed as a special case of the first part, with trivial character as a simple character, $F=F[\beta]$ and $\curlyb{\beta}=\curlyu$.
                                        
\subsection{Supercuspidal representations}
Let $(J,\lambda)$ be a simple type, with the simple stratum $[\curlyu,n,0,\beta]$ and $E=F[\beta]$.

\subsubsection{Maximal simple types} The following are equivalent:
\begin{itemize}
\item[(i)] $\eBoe{\beta}=1$
\item[(ii)] There exists an irreducible supercuspidal representation $\pi$ of $G$, such that $\Rest{G}{J}{\pi}$ contains $\lambda$
\item[(iii)]Any irreducible representation $\pi$ of $G$, such that $\Rest{G}{J}{\pi}$ contains $\lambda$, is supercuspidal.
\end{itemize}
Suppose these conditions hold, and let $\pi$ be an irreducible representation of $G$, which contains $\lambda$. Then an irreducible representation $\pi'$ will contain $\lambda$ if and only if $\pi'\cong\pi\otimes\chi\circ\det$, for some unramified quasicharacter $\chi$ of $F^{\times}$. We say that such $(J,\lambda)$ is a \textit{maximal} simple type, see \cite{bk}(6.2.3). 

We also note, that if $(J,\lambda)$ is a maximal simple type then
$$\IG{\lambda}{\lambda}{J}=\EJ,$$
see \cite{bk}(5.5.11) and (6.2.1).
\subsubsection{Structure of supercuspidal representations}

Let $\pi$ be an irreducible supercuspidal representation of $G$. There exists a simple type $(J,\lambda)$ in $G$, such that $\Rest{G}{J}{\pi}$ contains $\lambda$. Further, 
\begin{itemize}\item[(i)] the simple type $(J,\lambda)$ is uniquely determined up to $G$-conjugacy.
\item[(ii)]if $(J,\lambda)$ is given by a simple stratum $[\curlyu,n,0,\beta]$ with $E=F[\beta]$, there is a uniquely determined representation $\Lambda$ of $\EJ$, such that $\Rest{G}{J}{\Lambda}\cong\lambda$ and $\pi\cong\cIndu{\EJ}{G}{\Lambda}$.
\item[(iii)] If $(J,\lambda)=(\uz,\sigma)$, i.e., $J=\uz$ for some principal $\oF$ order $\curlyu$ and $\lambda$ is trivial on $\uu{1}$, then there exists a uniquely determined representation $\Lambda$ of $F^{\times}\uz$, such that $\Rest{G}{F^{\times}\uz}{\Lambda}\cong\lambda$ and $\pi\cong\cIndu{F^{\times}\uz}{G}{\Lambda}$.
\end{itemize}  
See \cite{bk}(8.4.1). Here $\cInd$ denotes compact induction, which is described in detail in  \cite{bu2}.

\subsubsection{Split types} 

A split type is a pair $(K',\vartheta)$, where $K'$ is a compact open subgroup of $G$, and $\vartheta$ is an irreducible representation $K'$. There are four flavours of split types, and we will define the three, that we require in the course of the paper. We will need to use the following result: 

Let $\pi'$ be a smooth irreducible representation of $G$. If $\Rest{G}{K'}{\pi'}$ contains $\vartheta$, then the Jacquet module $\pi'_{U}$ is nontrivial for some unipotent radical of a proper parabolic subgroup of $G$, see \cite{bk}(8.2.5) and (8.3.3). 

In particular, a supercuspidal representation cannot contain a split type.  

\subsection{The Setup}\label{setup}

In this paper every hereditary $\oF$ order $\curlyu$ in $A$ will come with a lattice chain:
$$\mathcal L:\ldots L_{i+1}\subset L_i\subset L_{i-1}\ldots$$
such that $\curlyu=\End_{\oF}(\mathcal L)$. The lattice chain $\mathcal L$ will come with a basis $v_1,\ldots,v_N$ of $V$, with respect to which $\curlyu$ is identified with the ring of block upper triangular matrices modulo $\pF$. If $\curlyu$ is principal, then:$$L_0=\oF v_1+\ldots+\oF v_N$$ $$L_i=\oF v_1+\ldots+\oF v_{\frac{N}{e}(e-i)}+\pF v_{\frac{N}{e}(e-i)+1}+\ldots+\pF v_N$$ for $0<i<e=\eUof$. We define a useful element $\Pi$ on the basis of $V$.      
$$\begin{aligned}
 \Pi: v_i&\mapsto\pi_F v_{{\frac{N}{e}(e-1)+i}}\mbox{, for }1\le i\le \frac{N}{e}\\
      v_j&\mapsto v_{j-{\frac{N}{e}}}\mbox{, for }\frac{N}{e}+1\le j\le N.
\end{aligned} $$
By inspecting how $\Pi$ acts on the lattice chain $\mathcal L$, we see that $\Pi\in\kU$ and $\valu(\Pi)=1$. Hence we have a short exact sequence:
\[\begin{CD}{1}@>>>{\uz}@>>>{\kU}@>{\valu}>>{\ZZ}@>>>{0}  \end{CD} \]
We will always denote $$K=\Aut_{\oF}(L_0)$$
So $\uz$ is always a subgroup of $K$ and with respect to our basis $K$ is identified with $\GL_N(\oF)$.

Throughout the paper we fix a supercuspidal representation $\pi$ of $G$. Let $(J,\lambda)$ be a simple type occurring in $\pi$, with a simple stratum $[\curlyu,n,0,\beta]$, $E=F[\beta]$. We define $$\rho=\Indu{J}{\uz}{\lambda}$$
Since $\IG{\lambda}{\lambda}{J}=E^{\times}J$, and $E^{\times}J\cap\uz=J$, as $J$ is the unique maximal compact open subgroup of $E^{\times}J$, the representation $\rho$ is irreducible. It is worth writing out the details, since we will use this kind of argument a lot:
$$\brac{\Indu{J}{\uz}{\lambda}}{\Indu{J}{\uz}{\lambda}}{\uz}=\brac{\lambda}{\Indu{J}{\uz}{\lambda}}{J}=\ldots$$
$$\ldots=\sum_{u\in J \backslash \uz/J}\brac{\lambda}{\Indu{J\cap J^u}{J}{\lambda^u}}{J}=\sum_{u\in J \backslash \uz/J}\brac{\lambda}{\lambda^u}{J\cap J^u}$$
the equalities above involve Frobenius reciprocity and Mackey's formula, hence $\brac{\rho}{\rho}{\uz}=1$. Since $[\curlyu,n,0,\beta]$ is simple, $\curlyu$ is principal and since $(J,\lambda)$ is contained in a supercuspidal representation, we have $\eBoe{\beta}=1$, hence $\valu(\pi_E)=1$. That implies $$\kU=E^{\times}\uz$$
Let $\Lambda$ be the unique extension of $\lambda$ to $\EJ$ such that $\pi\cong\cIndu{\EJ}{G}{\Lambda}$. We define $$\tilde{\rho}=\Indu{\EJ}{\kU}{\Lambda}$$
Just by transitivity of induction, we have $$\pi\cong\cIndu{\kU}{G}{\tilde{\rho}}$$ 
Since $\pi$ is irreducible, $\tilde{\rho}$ is also irreducible and $$\Rest{G}{\uz}{\tilde{\rho}}\cong{\rho}$$ since in this case there is only one double coset in Mackey's formula.

Now we forget all about our original $(J,\lambda)$. There are two justifications for this. The vague is: since $\pi\cong\cIndu{\kU}{G}{\tilde{\rho}}$, we do not lose any information. The rigorous one is: suppose $(J_1,\lambda_1)$ is another simple type contained in $\pi$, with a simple stratum $[\curlyu,n_1,0,\beta_1]$, then there exists $x\in\uz$, such that $(J,\lambda)=(J_1^x,\lambda_1^x)$, and hence $\rho\cong\Indu{J_1}{\uz}{\lambda_1}$ and $\tilde{\rho}\cong\Indu{E_1^{\times}J_1}{\kU}{\Lambda_1}$. To see that, one needs to go through the proof of \cite{bk}(5.7.1), which says "intertwining implies conjugacy", and in the last step use that $\eBoe{\beta}=1$. As this does not affect us, we will not provide the details.

Why do we prefer working with $\rho$, rather than with a simple type? We are interested in the irreducible summands of $\Rest{G}{K}{\pi}$, and
 $$\Rest{G}{K}{\pi}\cong \mackey{g}{K}{G}{\kU}{{\tilde{\rho}}}$$
 Since, $\uz$ is the unique maximal compact open subgroup of $\kU$
$$\Rest{G}{K}{\pi}\cong \Mplus{g}{K}{G}{{\kU}}\IndResg{g}{K}{G}{\uz}{\rho}$$
In order to acquire some information about the irreducible summands, we will need to choose some nice representative of a double coset $Kg\kU$. Now $\kU$ is reasonable to work with, since we can identify $\curlyu$ with block upper triangular matrices modulo $\pF$, and $\valu(\Pi)=1$, so $\Pi$ and $\uz$ will generate $\kU$. On the other hand, it would be a lot harder to work with double cosets $Kg\EJ$, just ask: for what matrices $\beta$ is $[\curlyu,n,0,\beta]$ a simple stratum, if $\curlyu$ is identified with the block upper triangular matrices modulo $\pF$?

Since, by the argument above, $\rho$ determines the restriction of $\pi$ to any compact open subgroup of $G$, we will often omit $\pi$ from the statements of propositions, and will work with $\rho$ instead.

\section{Existence}\label{existence}

\begin{prop}\label{exist} Let $\tau=\Indu{\uz}{K}{\rho}$, then $\tau$ is a type for $\Ip$. Moreover, $\tau$ occurs in $\Rest{G}{K}{\pi}$ with multiplicity one.
\end{prop}
\begin{proof}Let $(J,\lambda)$ be a simple type, with a simple stratum $[\curlyu,n,0,\beta]$, such that $\rho\cong\Indu{J}{\uz}{\lambda}$. Let $E=F[\beta]$. By \cite{bk}(5.5.11) coupled with \cite{bk}(6.2.1) we know $$\IG{\lambda}{\lambda}{J}=E^{\times}J\mbox{, and }E^{\times}J\cap K=J$$ since $J$ is the unique maximal open compact subgroup of $\EJ$, hence 
$\tau\cong\Indu{J}{K}{\lambda}$ is irreducible.
$$\Rest{G}{K}{\pi}\cong \Mplus{g}{K}{G}{{\kU}}\IndResg{g}{K}{G}{\uz}{\rho}$$
So $\tau$ occurs in $\Rest{G}{K}{\pi}$, and corresponds to the double coset $K\kU$. Restriction to $K$ forgets tensoring with unramified quasicharacters, so if $\pi'\in\Ip$, then $\Rest{G}{K}{\pi'}$ contains $\tau$.  

Suppose $\pi'$ contains $\tau$, then restriction to $J$ will contain $\lambda$ and by \cite{bk}(6.2.3) $\pi'\cong\pi \otimes \chi \circ \det_A$, for some unramified quasicharacter $\chi$ of $F^{\times}$.

If $\tau$ was contained in $\pi$  more than once, then by restricting to $J$, we would get that $\lambda$ was contained in $\pi$ more than once.  
$$\Rest{G}{J}{\pi}\cong \Mplus{h}{J}{G}{{\EJ}}\IndResg{h}{J}{G}{J}{\lambda}$$
Since $\IG{\lambda}{\lambda}{J}=\EJ$, $\lambda$ will only be a summand of the representation coming from the double coset $J.1.\EJ$, which is isomorphic to $\lambda$. So $\lambda$ occurs with multiplicity one, and hence $\tau$ occurs with multiplicity one.
\end{proof}

\section{Key}\label{important}

Suppose $\tau$ is any  irreducible representation of $K$ occurring in $\pi$.  
Since$$\Rest{G}{K}{\pi}\cong \Mplus{g}{K}{G}{{\kU}}\IndResg{g}{K}{G}{\uz}{\rho}$$
then for some $g\in Kg\kU$ we have $\brac{\tau}{\Indu{{K\cap \cgate{g}{\uz}}}{K}{\rho^g}}{K}\neq 0$.

Let $(J,\lambda)$ be a simple type, with the simple stratum $[\curlyu,n,0,\beta]$, such that $\rho\cong\Indu{J}{\uz}{\lambda}$.

 $$\Rest{\uz}{{\tcap{{g^{-1}}}{\uz}{K}}}{\rho}\cong \Mplus{u}{{\tcap{{g^{-1}}}{\uz}{K}}}{\uz}{J}\Indu{{\tcap{{g^{-1}}}{\cgate{u}{J}}{K}}}{{\tcap{{g^{-1}}}{{\uz}}{K}}}{\lambda^u}$$   
Hence $$\Indu{{\tcap{g}{K}{\uz}}}{K}{\rho^g}\cong\Mplus{u}{{\tcap{{g^{-1}}}{\uz}{K}}}{\uz}{J}\Indu{{\tcap{{ug}}{K}{J}}}{K}{\lambda^{ug}}$$
We have the freedom to replace $(J,\lambda)$, with $(J^u,\lambda^u)$, for any $u\in \uz$ and $\tau$ is a summand of at least one of the representations on the right. So we have the following proposition:

\begin{prop}\label{key} 
Suppose $\tau$ is an irreducible representation of $K$ and $g$ is a fixed representative of $Kg\kU$, such that $\brac{\tau}{\Indu{{K\cap \cgate{g}{\uz}}}{K}{\rho^g}}{K}\neq 0$, then there exists a simple type $(J,\lambda)$, with the simple stratum $[\curlyu,n,0,\beta]$, such that $$\rho\cong\Indu{J}{\uz}{\lambda}\mbox{ and }\brac{\tau}{\Indu{{K\cap \cgate{g}{J}}}{K}{\lambda^g}}{K}\neq 0$$
Moreover, suppose that for every irreducible summand $\xi$ of $\Rest{J}{{\tcap{{g^{-1}}}{J}{K}}}{\lambda}$ there exists a smooth irreducible representation $\lambda'$ of $J$, such that $$\brac{\xi}{\lambda'}{ \tcap{{g^{-1}}}{J}{K}}\neq 0\mbox{ and }\IG{\lambda'}{\lambda}{J}=\emptyset$$ or a smooth irreducible representation $\theta'$ of $H^1=\HibU{1}$, such that $$\Rest{H^1}{{\tcap{{g^{-1}}}{H^1}{K}}}{\theta'}=\Rest{H^1}{{\tcap{{g^{-1}}}{H^1}{K}}}{\theta}\mbox{ and }\IG{\theta'}{\theta}{H^1}=\emptyset$$ where $\theta\in \CUm{0}{\beta}$ and $\brac{\theta}{\lambda}{H^1}\neq 0$, then $\tau$ cannot be a type.
\end{prop}
\begin{proof} The first part of proposition is immediate from above. Suppose $\tau$ is a type. Let $\tilde{\tau}$ be an extension  of $\tau$ to $\FK$, such that $\Rest{G}{{\FK}}{\pi}$ contains 
$\tilde{\tau}$, then according to \cite{bk1}(5.2), $$\cIndu{{\FK}}{G}{\tilde{\tau}} \cong \coprod \pi \otimes \chi_i\circ \det$$ where $\chi_i$ are finitely many unramified quasicharacters of $F^{\times}$. So $\Rest{G}{J}{\pi}$ will contain all irreducible representations occurring in 
$$\cRest{J}{{\cIndu{{\FK}}{G}{\tilde{\tau}}}}\cong \Mplus{h}{J}{G}{{\FK}}\IndResg{h}{J}{G}{{K}}{\tau}$$ 
and $\Rest{G}{{H^1}}{\pi}$ will contain all irreducible representations occurring in 
$$\cRest{{H^1}}{{\cIndu{{\FK}}{G}{\tilde{\tau}}}}\cong \Mplus{h}{H^1}{G}{{\FK}}\IndResg{h}{{H^1}}{G}{{K}}{\tau}$$ 
Now $$\brac{\tau}{\lambda^g}{{K\cap \cgate{g}{J}}}=\brac{\tau}{\Indu{{K\cap \cgate{g}{J}}}{K}{\lambda^g}}{K}\neq 0$$ so $$\brac{\lambda}{\tau^{g^{-1}}}{{J\cap \cgate{{g^{-1}}}{K}}}\neq 0$$
Let $\xi$ be an irreducible representation of $\tcap{{g^{-1}}}{J}{K}$, such that
  $$\brac{\lambda}{\xi}{\tcap{{g^{-1}}}{J}{K}}\neq 0 \mbox{ and } \brac{\xi}{\tau^{g^{-1}}}{\tcap{{g^{-1}}}{J}{K}}\neq 0$$
By assumption there exists $\lambda'$, such that $$\brac{\lambda'}{\xi}{\tcap{{g^{-1}}}{J}{K}}\neq 0\mbox{ and } \IG{\lambda'}{\lambda}{J}=\emptyset$$ 
So $$\brac{\lambda'}{\Indu{{J\cap \cgate{{g^{-1}}}{K}}}{J}{\tau^{g^{-1}}}}{J}=\brac{\lambda'}{\tau^{g^{-1}}}{{J\cap \cgate{{g^{-1}}}{K}}}\neq 0$$
That implies $\lambda'$ occurs in $\Rest{G}{J}{\pi}$. On the other hand we know that for some extension $\Lambda$ of $\lambda$ to $\EJ$, we have $\pi\cong\cIndu{{\EJ}}{G}{\Lambda}$, so 
$$\Rest{G}{J}{\pi}\cong \Mplus{h}{J}{G}{{\EJ}}\IndResg{h}{J}{G}{J}{\lambda}$$
That implies that $\lambda$ and $\lambda'$ must intertwine in $G$, which is a contradiction.

We deal similarly with $\theta'$. If $\lambda\cong\kappa\otimes\sigma$, then by unravelling all the definitions we have 
$$\Rest{J}{{H^1}}{\lambda}\cong(\dim \sigma \dim \kappa)\theta$$
Hence $$\brac{\theta}{\tau^{g^{-1}}}{{H^1\cap \cgate{{g^{-1}}}{K}}}\neq 0$$
Then by the same argument as above, we show that $\Rest{H^1}{{\tcap{{g^{-1}}}{H^1}{K}}}{\theta'}=\Rest{H^1}{{\tcap{{g^{-1}}}{H^1}{K}}}{\theta}$ implies that $\theta'$ occurs in $\Rest{G}{{H^1}}{\pi}$. $$\Rest{G}{{H^1}}{\pi}\cong \Mplus{h}{{H^1}}{G}{{\EJ}}\IndResg{h}{{H^1}}{G}{J}{\lambda}$$ 
Hence $\brac{\theta'}{\lambda^h}{{H^1\cap \cgate{h}{J}}}\neq 0$, for some $h\in G$, so $\brac{\theta'}{\lambda^h}{{H^1\cap \cgate{h}{{H^1}}}}\neq 0$, which implies $h\in\IG{\theta'}{\theta}{H^1}$. And we obtain a contradiction.  
\end{proof}
\begin{remar}The conditions of the Proposition above might seem a little strange. So we give the following example. Suppose $\lambda$, $\lambda'$ as above, and assume further, that $(J,\lambda')$ is a simple type, with a simple stratum $[\curlyu,n',0,\beta']$, occurring in a supercuspidal representation $\pi'$. Now the condition on irreducible summands translates into $$\brac{\tau}{\Indu{K\cap J^{g}}{K}{\lambda{'}^{g}}}{K}\neq 0\mbox{ and }\brac{\tau}{\Indu{K\cap J^{g}}{K}{\lambda^{g}}}{K}\neq 0$$
Hence $$\brac{\tau}{\pi'}{K}\neq 0 \mbox{ and }\brac{\tau}{\pi}{K}\neq 0$$
And since $\lambda$ and $\lambda'$ do not intertwine, we have $$\pi'\not\cong\pi\otimes\chi\circ\det$$
where $\chi$ is an unramified quasicharacter of $F^{\times}$. So $\tau$ cannot be a type. If $N=2$, this situation arises in \cite{ap}\S A.3.7 and \S A.3.10.
\end{remar}
The rest of the paper is concerned with picking a nice representative from a double coset $Kg\kU$, and constructing $\lambda'$ and $\theta'$, when $Kg\kU \neq K\kU$.

\section{Representatives of  $Kg\kU$}\label{represent}
We will need to identify $\curlyu$ with the ring of block upper triangular matrices modulo $\pF$ in order to do explicit calculations. For that purpose we introduce the following notation.
\begin{nota} We will write $\MM(m,\oF)$ for the ring of $m\times m$ matrices with coefficients in $\oF$. We will also write $\MM(m,\pF^i)=\pi_F^i\MM(m,\oF)$.
\end{nota}
\begin{prop}\label{cosets}Let $\curlyu$ be a principal hereditary $\oF$ order in $A$, and let $e=\eUof$. Suppose a double coset $Kg\kU\neq K\kU$, then there exists a representative $g$ of $Kg\kU$, such that one of the following holds:
\begin{enumerate}
\item[(A)]  The map $\uz\cap \cgate{{g^{-1}}}{K}\rightarrow\uz/\uu{1}$ is not surjective. Moreover, for some  index $j$, $0\leq\ j <e$, the image of $\uz\cap \cgate{{g^{-1}}}{K}$ in $\Aut_{\kF}(L_j/L_{j+1})$, via the map:
\[
\begin{CD}
 {\uz\cap \cgate{{g^{-1}}}{K}}@>>> {\uz} @>>> {\Aut_{\kF}(L_j/L_{j+1})}\\
\end{CD} 
\]
is a proper parabolic subgroup of $\Aut_{\kF}(L_j/L_{j+1})$. 
\item[(B)] The map $\uz\cap \cgate{{g^{-1}}}{K}\rightarrow\uz/\uu{1}$ is  surjective, and $$(h-1).L_{e-1}\subseteq L_{e+1}\mbox{, }\every h \in\uu{1}\cap \cgate{{g^{-1}}}{K} $$
\end{enumerate}
\end{prop}
\begin{proof}We identify $\curlyu$ with the ring of block upper triangular matrices modulo $\pF$, with respect to our basis $v_1,\ldots,v_N$ of $V$. Then $K$ is identified with $\GL_N(\oF)$. Having made these identifications, we prove the lemma below:
\begin{lem}There exists a representative $g\in Kg\kU$, such that $g$ is a diagonal matrix and the diagonal entries $=(\uni{1},\ldots,\uni{N})$, where $$\al{{i\frac{N}{e}+1}}\geq \ldots \geq \al{{(i+1)\frac{N}{e}}}\geq 0$$ for all $0\leq i < e$ and one of the following holds:
\begin{enumerate}
\item[(A)]  $\al{{j\frac{N}{e}+1}}\neq\al{{(j+1)\frac{N}{e}}}$, for some  $j$, $0\leq\ j <e$.
\item[(B)] $\al{{i\frac{N}{e}+1}}=\al{{(i+1)\frac{N}{e}}}$, for all  $i$, $0\leq\ i <e$, $\al{1}\geq 2$ and there exists an index $j$, such that $\al{k}>0$ if $k<j$ and 
$\al{k}=0$ if $k\geq j$, for all $1\leq k\leq N$.
\end{enumerate}
\end{lem}
\begin{proof} Let $\curlyu_m\subseteq\curlyu$ be the upper triangular matrices modulo $\pF$. Then the Iwahori decomposition tells us that $G$ is a disjoint union of double cosets $\mathbf U(\curlyu_m)w\mathbf U(\curlyu_m)$, for $w\in\tilde{W}=W_0\ltimes D$, where $W_0$ is the group of permutation matrices and $D$ is is the group of diagonal matrices, whose eigenvalues are powers of $\pi_F$.
We have $W_0\le K$ and using permutation matrices in $K$ and $\uz$, we can choose $g$ to be diagonal with diagonal entries $=(\uni{1},\ldots,\uni{N})$, where $\al{{i\frac{N}{e}+1}}\geq \ldots \geq \al{{(i+1)\frac{N}{e}}}$, for all $0\leq i <e$. By multiplying $g$ by an element of $F^{\times}$, we can ensure that $\al{i}\ge 0$ for all $1\le i\le N$ and at least one of them is equal to $0$. If $(A)$ is true, then we are done. 

Otherwise, let $\Pi\in \kU$ be the element defined in Section \ref{setup}, and
let $t$ be the following  permutation matrix: 
$$\begin{aligned}
 t: v_{{\frac{N}{e}(e-1)+i}}& \mapsto v_i \mbox{, for }1\le i\le \frac{N}{e}\\
      v_{j-{\frac{N}{e}}}& \mapsto v_j\mbox{, for }\frac{N}{e}+1\le j\le N.
\end{aligned} $$
We write $(\al{1},\ldots,\al{e})$ for the diagonal matrix $(\uni{1}I,\ldots,\uni{i}I,\ldots,\uni{e}I)$, and $I$ is the  $\frac{N}{e}\times\frac{N}{e}$ identity matrix. Let $\oplus$ be the map $\oplus:g\mapsto tg\Pi$, then $$\oplus:(\al{1},\ldots,\al{e})\mapsto(\al{e}+1,\al{1},\ldots,\al{{e-1}})$$  Similarly, let $\ominus$ be the map $\ominus:g\mapsto t^{-1}g\Pi^{-1}$, then $$\ominus:(\al{1},\ldots,\al{e})\mapsto(\al{2},\ldots,\al{e},\al{1}-1)$$    Let $j$ be the smallest index such that $\al{j}=0$, then we replace $g$ with $\oplus(g)$ $e-j$ times. If all $\al{i}\le 1$, then by replacing $g$ with $\ominus(g)$ $e-1$ times, we would obtain the identity matrix. The double coset $Kg\kU\neq K\kU$, so let $k$ be the smallest index such that $\al{k}\ge 2$, then by replacing $g$ with $\ominus(g)$ $k-1$ times we get $g$ of the required form. 
\end{proof}
If $g$ satisfies part (A) of the Lemma, then $g$ will also satisfy part (A) of the proposition. To see that it is enough to know what the  $\frac{N}{e}\times\frac{N}{e}$ blocks on the diagonal of a matrix in $K\cap \cgate{{g^{-1}}}{K}$ look like. Since $g$ is chosen nicely, it is enough to do the computation for $e=1$ and then reduce modulo $\pF$, which is easy.

If $g$ satisfies part (B) of the Lemma, then $g$ will also satisfy part (B) of the Proposition. Surjectivity follows by the same argument as in (A). For the second part we observe that every matrix $(A_{ij})\in \uu{1}\cap \cgate{{g^{-1}}}{K}$, where $A_{ij}\in\MMe{\oF}$, has $A_{e1}\in\MMe{\pF^2}$. It is enough to do the calculation for $K\cap \cgate{{g^{-1}}}{K}$ and since $g$ is chosen nicely, we may assume $e=N$, and then it is obvious. Since with respect to our fixed basis $v_1,\ldots,v_N$ of $V$:
$$L_{e-1}=\oF v_1+\ldots+\oF v_{\frac{N}{e}}+\pF v_{\frac{N}{e}+1}+\ldots+\pF v_N$$
$$L_e=\pF v_1+\ldots+\pF v_N$$
$$L_{e+1}=\pF v_1+\ldots+\pF v_{\frac{N}{e}(e-1)}+\pF^2 v_{\frac{N}{e}(e-1)+1}+\ldots+\pF^2 v_N$$ 
we have  $$(h-1).L_{e-1}\subseteq L_{e+1}\mbox{, }\every h \in\uu{1}\cap \cgate{{g^{-1}}}{K} $$
\end{proof}
\begin{defi} Let $\curlyu$ be a principal hereditary $\oF$ order in $A$ and $g$ be a representative of a double coset  $Kg\kU$. We say a pair $(g,Kg\kU)$ has \textbf{property (A)} (resp. \textbf{property (B)}) if $Kg\kU\neq K\kU$ and $g$ satisfies \ref{cosets} (A) (resp. \ref{cosets} (B)). 
\end{defi}

\begin{remar} If $\eUof=1$, then all the double cosets have property (A). If $\eUof=N$, then all the double cosets have property (B). In particular, when $N=2$, the two cases above correspond to the ones considered in \cite{ap}. 
\end{remar}

\section{Double cosets with property (A)}\label{a}
If $(g,Kg\kU)$ has property (A), then the map $\uz\cap \cgate{{g^{-1}}}{K}\rightarrow\uz/\uu{1}$ is not surjective.  If $N=2$, then $(g,Kg\kU)$ has property (A) if and only if $\eUof=1$, in which case $\uz\cap \cgate{{g^{-1}}}{K}$ coincides with the group considered in \cite{ap}\S A.3.7. For general $N$ we show that if $[\curlyu,n,0,\beta]$ is a simple stratum, $E=F[\beta]$ and $\eBoe{\beta}=1$, then the map $$J(\beta,\curlyu)\cap \cgate{{g^{-1}}}{K}\rightarrow J(\beta,\curlyu)/J^1(\beta, \curlyu)\cong \ubz{\beta}/\ub{1}{\beta}$$ is not surjective. So we look for $\lambda'$, satisfying the conditions of Proposition \ref{key}, of the form $\lambda'=\kappa\otimes\sigma'$, where $\sigma'$ is a lift to $J$ of an irreducible representation of $J/J^1\cong\ubz{\beta}/\ub{1}{\beta}$. Since $J^1\le\Ker\sigma$, the restriction $\Rest{J}{{\tcap{{g^{-1}}}{J}{K}}}{\sigma}$ depends only on the image of  $\tcap{{g^{-1}}}{J}{K}$ in $\ubz{\beta}/\ub{1}{\beta}$.
\begin{defi} Suppose that $(g,Kg\kU)$ has property (A), then we define $\Kg$ to be the  following subgroup of $\uz$: $$\Kg=(\uz\cap \cgate{{g^{-1}}}{K})\uu{1}.$$ 
\end{defi}
From the definition we get that $\Kg$ is a parahoric subgroup of $G$, contained in $\uz$. Since $ \uu{1}$ is a subgroup of $\Kg$, we have  $$\Kg\cap\JbU=(\Kg\cap\ubz{\beta})\JibU{1}.$$ 
For each $i$, $\ubz{\beta}$ acts on a $\kE$ vector space $L_i/L_{i+1}$ and $\uz$ acts on a $\kF$ vector space $L_i/L_{i+1}$. Let $j$ be the index, such that the image of $\Kg$ in $\Aut_{\kF}(L_j/L_{j+1})$ is a proper parabolic subgroup $P$. Let $H$ be the image of $\Kg\cap\ubz{\beta}$ in $\Aut_{\kE}(L_j/L_{j+1})$. We have the following diagram: 
\[
\begin{CD}
{\Kg\cap\ubz{\beta}} @>>> {\ubz{\beta}} @>>> {\Aut_{\kE}(L_j/L_{j+1})}\\
@VVV             @VVV       @VVV \\
{\Kg} @>>> {\uz} @>>> {\Aut_{\kF}(L_j/L_{j+1})}\\
\end{CD} 
\]
Horizontal arrows on the right are surjections, and all the other arrows are inclusions. Hence we get an injection $H\hookrightarrow P$. This injection will give us enough information about $H$, to handle $\Rest{G}{\Kg\cap\JbU}{\sigma}$.  We make the following definition.

\begin{defi}Let $\FF_q$ be a finite field, $H$ a subgroup $\GL_N(\FF_q)$, $\FF_{q^N}$ the unique extension of $\FF_q$ of degree $N$, and  $$S=\{h\in \GL_N(\FF_q):\chi_h(X)=f(X)^l, l\in\NN, f(X)\mbox{ irreducible over }\FF_q\}$$where $\chi_h(X)$ is a characteristic polynomial of $h$. Call $H$ \textbf{sufficiently small} if there exists a subfield $\FF$ of $\FF_{q^N}$, such that $[\FF_{q^N}:\FF]>1$ and  for all $h \in H\cap S$ the roots of $\chi_h(X)$ lie in $\FF$.
\end{defi}

\begin{lem}\label{parsmall} Let $P$ be a proper parabolic subgroup of $\GL_N(\FF_q)$, then $P$ is sufficiently small.
\end{lem}
\begin{proof}Without loss of generality we may assume that $P$ is a maximal proper parabolic subgroup. As conjugation does not change the characteristic polynomial, we may further assume that $P$ is a subgroup of  block upper triangular matrices  consisting of two blocks of size $a\times a$ and $b\times b$, where $a+b=N$. If $h\in P$, then  the characteristic polynomial $\chi_h(X)$ of $h$ can be written as a product $\chi_h(X)= f_1(X)f_2(X)$, where $\deg(f_1)=a$ and $\deg(f_2)=b$. Hence if $h\in P$ and $\chi_h(X)= f(X)^l$, where $f(X)$ is irreducible over $\FF_q$, then $\deg(f)$ divides $a$ and $\deg(f)$ divides $b$, so the roots of $f(X)$ lie in $\FF_{q^c}$, where $c=\gcd(a,b)$. As $c$ divides $N$ and $c < N$, we deduce that $P$ is sufficiently small. 
\end{proof}
\begin{remar}One might think that every sufficiently small subgroup is contained in a proper parabolic subgroup. The following example shows that this is not the case. Choose $a$, $N>1$, such that $\gcd(a,N)=1$, then $\GL_N(\FF_q)$ is a sufficiently small subgroup of $\GL_N(\FF_{q^a})$, with $\FF=\FF_{q^N}$. It cannot be contained in any proper parabolic subgroup of $\GL_N(\FF_{q^a})$, since $\FF_{q^{Na}}$ is the smallest extension of $\FF_{q^a}$ containing $\FF_{q^N}$. It is also not hard to construct an embedding $\iota:\GL_N(\FF_{q^a})\hookrightarrow \GL_{Na}(\FF_{q})$, such that $\iota(\GL_N(\FF_q))$ is contained in a proper parabolic subgroup of $GL_{Na}(\FF_{q})$, which is the case considered in the Lemma below. 
\end{remar}
\begin{lem}\label{minpol}Let $W$ be an $\FF_{q^a}$ vector space of dimension $N$, which we also consider as an $\FF_q$ vector space. So we get an embedding of algebras $$\iota:\End_{\FF_{q^a}}(W)\hookrightarrow\End_{\FF_q}(W)$$ Let $H$ be a subgroup of $\Aut_{\FF_{q^a}}(W)$, such that $\iota(H)$ is contained in a proper parabolic subgroup $P$ of $\Aut_{\FF_q}(W)$, then $H$ is a sufficiently small subgroup of $\Aut_{\FF_{q^a}}(W)$. 
\end{lem}

\begin{proof} Let $h\in H$ and suppose that the characteristic polynomial $\chi_h(X)$ of $h$ in $\End_{\FF_{q^a}}(W)$ is a power of a polynomial $f(X)$, which is irreducible over $\FF_{q^a}$. Let $\FF_{q^b}$ be the field  generated by the coefficients of $f(X)$ over $\FF_q$. Define $\tilde{f}(X)$ to be $$\tilde{f}(X)=\prod_{\xi\in \Gal(\FF_{q^b}/\FF_q)}f(X)^{\xi}$$ where $\Gal(\FF_{q^b}/\FF_q)$ acts on the coefficients of $f(X)$. Then $\tilde{f}(X)\in\FF_q[X]$ and $\tilde{f}(X)$ is irreducible over $\FF_q$. We claim that the characteristic polynomial of $\iota(h)$ in $\End_{\FF_q}(W)$ is a power of $\tilde{f}(X)$. To see that it is enough to show that the minimal polynomial of $\iota(h)$ divides some power of $\tilde{f}(X)$, but $\tilde{f}(\iota(h))=\iota(\tilde{f}(h))$, and now the claim is obvious.

We apply Lemma \ref{parsmall} to $P$ and  $\Aut_{\FF_q}(W)$ and  hence we get  a subfield $\FF$ of $\FF_{q^{aN}}$, such that $[\FF_{q^{aN}}:\FF]>1$, and roots of $\tilde{f}(X)$ and hence roots of $f(X)$, lie in $\FF$. The field $\FF$ does not depend on the choice of $h$, so $H$ is a sufficiently small subgroup of $\Aut_{\FF_{q^a}}(W)$.      

\end{proof}

\begin{cor}\label{bananas} Suppose $(g,Kg\kU)$ has property (A). Let $[\curlyu,n,0,\beta]$  be a simple stratum, such that $\eBoe{\beta}=1$ and let $E=F[\beta]$. Let $H$ be the image of $\Kg\cap\ubz{\beta}$ in $\ubz{\beta}/\ub{1}{\beta}$. Then $H$ is a sufficiently small subgroup of $\ubz{\beta}/\ub{1}{\beta}$ and if $E$ is 
a totally ramified extension of $F$, then $H$ is contained in a proper parabolic subgroup of $\ubz{\beta}/\ub{1}{\beta}$. Moreover  $H\neq \ubz{\beta}/\ub{1}{\beta}$.
\end{cor}
\begin{proof} As $(g,Kg\kU)$ has property (A), we can find a lattice $L_j$ in the lattice chain defining $\curlyu$, such that the image of $\Kg$ in 
 $\Aut_{\kF}(L_j/L_{j+1})$ is a proper parabolic subgroup $P$. As $\eBoe{\beta}=1$, we get that $L_{j+1}=\pi_E L_j$ and hence $\ubz{\beta}/\ub{1}{\beta}\cong\Aut_{\kE}(L_j/L_{j+1})$, so Lemma \ref{minpol} applied to $L_j/L_{j+1}$, implies that $H$ is a sufficiently small subgroup of  
$\ubz{\beta}/\ub{1}{\beta}$.

If $E$ is a totally ramified extension of $F$, then $\kE=\kF$ and $H\le P$, so $H$ is contained in a proper parabolic subgroup. Note that, in this case $\dim_{\kE}(L_j/L_{j+1})=1$ would imply $\eUof=N$, and hence $(g,Kg\kU)$ has property (B).

We can pick a polynomial $f(X)\in\kE[X]$ of degree $\dim_{\kE}(L_j/L_{j+1})$,  which is irreducible over $\kE$ and $f(X)\neq X$. From linear algebra we know that there exists some $h\in\Aut_{\kE}(L_j/L_{j+1})$, such that the characteristic polynomial of $h$ equals to $f(X)$. That implies  $H\neq \ubz{\beta}/\ub{1}{\beta}$, as $H$ is sufficiently small.  
\end{proof}

\begin{lem}\label{zsigmondy} For all integers $q>1$ and $N>1$ there exists a prime $r$ such that $r$ divides $q^N-1$, but $r$ does not divide $q^m-1$, for all $0<m<N$, except when $q=2^i-1$ and $N=2$ or $q=2$ and $N=6$.
\end{lem}
\begin{proof}This result is known as Zsigmondy's theorem. We refer the reader to \cite{zs}.
\end{proof}

\begin{prop}\label{cuspidals}
Let $\sigma$ be a cuspidal irreducible  representation of $\GL_N(\FF_{q})$ affording a character $\XX$, where  $\FF_{q}$ is a finite field with  $q$
 elements and $p$ is the characteristic of  $\FF_{q}$.

 Suppose $H$ is a sufficiently small subgroup of $\GL_N(\FF_{q})$, and if $q=2$ or $q=3$, we further assume that $H$ is contained in a proper parabolic 
subgroup of $\GL_N(\FF_{q})$. Then for every irreducible representation $\xi$ of $H$, such that $\brac{\xi}{\sigma}{H}\neq 0$, 
there exists an irreducible representation $\sigma'$ of  $\GL_N(\FF_{q})$, such that $\sigma\not\cong\sigma'$ and $\brac{\xi}{\sigma'}{H}\neq 0$.

 Moreover, in all, except finitely many, cases we may choose $\sigma'$ to be a cuspidal representation,
 such that $\Rest{G}{H}{\sigma}\cong\Rest{G}{H}{\sigma'}$.
\end{prop}
\begin{remar} So $H$ is small enough to not distinguish between two different cuspidal representations.
\end{remar} 
\begin{proof} We denote $\GL_N(\FF_{q})$ by $\Gamma$ and let $$S=\{h\in \Gamma:\chi_h(X)=f(X)^l, l\in\NN, f(X)\mbox{ irreducible over }\FF_q\}$$ where $\chi_h(X)$ is a characteristic polynomial of $h$. Let $\FF_{q^N}$ be an extension of $\FF_q$ of degree $N$. As $H$ is a sufficiently small subgroup there exists a subfield $\FF$ of $\FF_{q^N}$, such that $[\FF_{q^N}:\FF]>1$ and for every $h\in H\cap S$ the roots of the characteristic polynomial of $h$ lie in $\FF$. First of all we get rid of some easy cases:

If $N=1$, then $\sigma$ is a one dimensional representation. Let $\Psi$ be a lift to $\FF_q^{\times}$ of some non-trivial linear character of $\FF_q^{\times}/\FF^{\times}$, then $\sigma'=\sigma\otimes\Psi$, satisfies the conditions of the proposition.

So we may assume that $N\ge 2$. The proposition is false if and only if we can find an irreducible representation $\xi$ of $H$ such that $\Indu{H}{\Gamma}{\xi}\cong \sigma \oplus \ldots \oplus\sigma$. Suppose $H$ is contained in some proper parabolic subgroup $P$. Let $U$ be the unipotent radical of the parabolic subgroup opposite to $P$. Then $U\cap P=1$ and hence $U\cap H=1$, so $\brac{\Eins}{\Indu{H}{\Gamma}{\xi}}{U}\neq 0$, which implies that $\brac{\Eins}{\sigma}{U}\neq 0$, but $\sigma$ is a cuspidal representation, so we obtain a contradiction.

In general we use character theory. The characters of the irreducible representations of $\Gamma$ were first described in \cite{green}, but \cite{mac} is also very useful. The conjugacy classes of $\Gamma$ are in one-to-one correspondence with isomorphism classes of $\FF_q[X]$ modules $W$ , such that $\dim_{\FF_q} W=N$ and $X.w=0$ implies that $w=0$, see \cite{mac}\S IV.2. Each $h\in \Gamma$ acts naturally on $\FF_q^N$, and hence defines an $\FF_q[X]$ module structure on $\FF_q^N$, such that $X.w=hw$, for all $w\in \FF_q^N$. We denote this module by $W_h$. Clearly, two elements $h_1$, $h_2\in\Gamma$ are conjugate if and only if  $W_{h_1}\cong W_{h_2}$. We may therefore write $W_c$ instead of $W_h$, where $c$ is the conjugacy class of $h$ in $\Gamma$. In our case we are only interested in those conjugacy classes $c$, where the characteristic polynomial $\chi_h(X)$ of any $h\in c$ is a power of a polynomial $f(X)$, which is irreducible over $\FF_q$. Since $\FF_q[X]$ is a principal ideal domain, and $\chi_h(X)$ will anniolate  $W_c$, we have
 $$W_c\cong\bigoplus_{i=1}^k \FF_q[X]/(f)^{\mu_i(c)}$$  That defines a partition $\mu(c)=(\mu_1(c), \mu_2(c),\ldots,\mu_k(c))$ of $\frac{N}{d}$, where $d$ is the degree of $f(X)$.

We are now ready to describe the characters of cuspidal representations of $\Gamma$, as given in \cite{gel1} and \cite{gel2}. Let $\Psi: \FF_{q^N}^{\times}\rightarrow \CC^{\times}$, be an abelian character, such that
 $\Psi^{{q^m-1}}\neq 1$ for all $m$ dividing but not equal to $N$, then the following class function $\XX_{\Psi}$ is a character of a cuspidal representation of $\Gamma$.
 \[
\XX_{\Psi}(h)=\left\{ \begin{array}
              {l@{\quad:\quad}l}
              0 & h\not\in S \\(-1)^{k+N}\varphi_k(q^d)(\Psi(\alpha^q)+\ldots+\Psi(\alpha^{q^d})) & h\in S
              \end{array} \right. \] 
where $\varphi_k(X)=(X-1)\ldots(X^{k-1}-1)$, $\varphi_1(X)=1$ and if $\chi_h(X)= f(X)^l$, $f(X)$ irreducible over $\FF_q$, then 
$d$ is the degree of $f$, $\alpha$ is a root of $f$ and $k$ is the number of parts in the partition $\mu(c)$, given by the conjugacy class $c$ of $h$. 
Conversely, any character $\XX$ of a cuspidal representation of $\GL_N(\FF_{q})$ arises in this way.  Moreover, if $\Theta$ is another abelian character $\Theta: \FF_{q^N}^{\times}\rightarrow \CC^{\times}$, such that $\Theta^{{q^m-1}}\neq 1$ for all $m$ dividing but not equal to $N$, then  $\XX_\Psi=\XX_\Theta$ if and only if $\Psi=\Theta^{q^m}$, for some $m\geq 0$.

Let $\Psi$ be an abelian character, such that $\XX=\XX_{\Psi}$ and suppose there exists an abelian character $\Theta:\FF_{q^N}^{\times}\rightarrow \CC^{\times}$, such that $\Rest{G}{\FF^{\times}}{\Psi}=\Rest{G}{\FF^{\times}}{\Theta}$, $\Theta^{q^m-1}\neq1$, for all $m$ dividing but not equal to $N$, and $\Theta\neq\Psi^{q^m}$, for $m\ge 0$, then we take $\sigma'$ to be a cuspidal representation of $\Gamma$ corresponding to the character $\XX_{\Theta}$. Since $\Theta\neq\Psi^{q^m}$, the representations $\sigma$ and $\sigma'$ are not isomorphic, and since $\XX_{\Psi}(h)=\XX_{\Theta}(h)$ for all $h\in H$, $\Rest{G}{H}{\sigma}$ is isomorphic to $\Rest{G}{H}{\sigma'}$.

In most cases, we can show such $\Theta$ exists, by counting characters with the  desired properties. The argument below was shown to me by S.\,D.\,Cohen. Let $a=[\FF_q:\FF_p]$ and $b=[\FF:\FF_p]$. By Lemma \ref{zsigmondy}, there will exist a prime $r$, such that $r$ divides $p^{aN}-1=q^N-1$, but $r$ does not divide $p^m-1$, for all $0<m<N$, unless $aN=2$ and $p=2^i-1$, or $aN=6$ and $p=2$. 

Suppose we are not in one of these exceptional cases. If $\Theta$ is not an $r$-th power, then $r$ divides the order of $\Theta$, and hence $\Theta^{p^m-1}\neq 1$, for all $0<m<aN$. In particular, $\Theta^{q^m-1}\neq 1$, for all $m$ dividing, but not equal to $N$. Since $\FF_{q^N}^{\times}$ is cyclic and $r$ does not divide $|\FF^{\times}|$, every abelian  character of $\FF^{\times}$ is a restriction of an abelian character of $\FF_{q^N}^{\times}$, which is an $r$-th power. Hence there will be   
$(1-\frac{1}{r})\frac{q^N-1}{p^b-1}$ characters $\Theta$, such that $\Theta$ is not an $r$-th power and $\Rest{G}{\FF^{\times}}{\Psi}=\Rest{G}{\FF^{\times}}{\Theta}$. In order to avoid $\Psi^{q^m}= \Theta$ for some $m\ge 0$, we need the following inequality to hold: $$(1-\frac{1}{r})\frac{q^N-1}{p^b-1}>N$$ Since $\FF$ is a proper subfield of $\FF_{q^N}$ we have $|\FF|\leq q^{\frac{N}{2}}$. So it is enough to prove $$(1-\frac{1}{r})(q^{\frac{N}{2}}+1)>N$$ That can be done using induction on $N$, when $q\geq 4$, then $r\geq 2$, when $q=3$ and $N\geq 2$, then $r \geq 5$. When $q=2$, $N\geq 3$, then $r \geq 5$, since  $3=2^2-1$.

We are left with the following cases: $\GL_2(\FF_2)$, $\GL_6(\FF_2)$, $\GL_3(\FF_4)$, $\GL_2(\FF_8)$ and $\GL_2(\FF_q)$, where $q$ is a prime and $q=2^i-1$. 

If $q=2$ or $q=3$, then by our assumption on $H$, it is contained in a proper parabolic subgroup and we have already dealt with this. 

If $N=2$ and $q=2^i-1$, where $q$ is a prime, $q>3$, then 
$\frac{q+1}{2}>3$. Therefore we may pick an abelian character $\Xi$ of $\FF_{q^2}^{\times}$, such that $\Rest{G}{\FF_q^{\times}}{\Xi}=1$ and $\Xi^2$ is not one of the following characters: $1$, $\Psi^{q-1}$ or $\Psi^{2(q-1)}$. Let $\Theta=\Psi\Xi$, then $\Theta^{q-1}=1$ implies $\Xi^2=\Psi^{q-1}$, $\Theta=\Psi$ implies $\Xi^2=1$ and $\Theta=\Psi^q$ implies $\Xi^2= \Psi^{2(q-1)}$.

If $N=3$ and $q=4$ or $N=2$ and $q=8$, let $\alpha$ be an element of order $63$ in $\FF_{64}^{\times}$ and let $c$ be the conjugacy class in $\Gamma$ corresponding $\FF_q[X]$ module $W=\FF_q[X]/(m_{\alpha})$, where $m_{\alpha}(X)$ is the minimal polynomial of $\alpha$ over $\FF_q$. If $\Indu{H}{\Gamma}{\xi}\cong \sigma \oplus \ldots \oplus\sigma$, for some representation $\xi$ of $H$, then $\XX(c)=0$, since $H$ does not meet $c$ in $\Gamma$, as $\alpha$ is not contained in any proper subfield of $\FF_{64}$. 

If $N=2$ and $q=8$, that implies $\Psi(\alpha)+\Psi(\alpha^8)=0$, so $2$ divides the order of $\Psi$ , which is impossible. 

If $N=3$ and $q=4$, then $\Psi(\alpha)+\Psi(\alpha^4)+\Psi(\alpha^{16})=0$, which implies $\Psi$ has order $9$, since $X^5+X+1=(X^3-X^2+1)(X^2+X+1)$ and $\Psi(\alpha)$ is a root of unity in $\CC$. If $\FF=\FF_8$ we may take $\Theta=\Psi^2$, since $\Rest{G}{\FF_8^{\times}}{\Theta}=\Rest{G}{\FF_8^{\times}}{\Psi}=1$ and $2\not\equiv 4^m \pmod{9}$. If $\FF=\FF_4$, then we may choose $\Theta=\Psi\Xi$, where $\Xi$ is any character of order $7$.    
\end{proof}
\begin{remar}When $\pi$ is a supercuspidal representation of $\GL_2(F)$, it is enough to prove the proposition above for $\GL_1(\FF_q)$, see \cite{ap} \S A.3.7.
\end{remar}
We recall the following definition.
\begin{defi}\label{splitone}\cite{bk}(8.1.2) A \textbf{split type of level $(0,0)$} is a pair $(K',\vartheta)$, given as follows:
\begin{itemize}
\item[(i)] $\curlyu$ is a hereditary $\oF$ order in $A$
\item[(ii)] $\uz/\uu{1}= \mathcal{G}_1\times\ldots\times\mathcal{G}_e$, where $\mathcal{G}_1\cong\GL_{n_i}(\kF)$
\item[(iii)] $K'=\uz$ and $\vartheta$ is the  inflation of an irreducible representation of $\uz/\uu{1}$ of the form $\xi_1\otimes\ldots\otimes\xi_e$, where $\xi_i$ is a cuspidal representation of $\mathcal{G}_i$, and $\xi_i\not\cong\xi_j$, for some $i\neq j$.
\end{itemize}
\end{defi} 
\begin{prop}\label{levelzero} Let $\pi$ be a supercuspidal representation of $G$ and $\curlyu$ a maximal hereditary $\oF$ order in $A$, i.e., $\eUof=1$. Suppose the space of vectors fixed by $\uu{1}$, $\pi^{\uu{1}}\neq 0$, then  $\pi^{\uu{1}}$ considered as a representation of $\uz$ is a lift of an irreducible cuspidal representation of $\uz/\uu{1}$.
\end{prop}
\begin{proof}If $\sigma$ is a representation of $\uz$, which is an irreducible summand of $\pi^{\uu{1}}$, then $\sigma$ is a lift of an irreducible representation of $\uz/\uu{1}$.
If $\sigma$ is not a lift of a cuspidal representation, then there exists a hereditary $\oF$ order $\curlyu'$ in $A$, such that $\curlyu'\subset\curlyu$, the image of 
$\mathbf U(\curlyu')$ in $\uz/\uu{1}$ is a proper parabolic subgroup
and $\Rest{\uz}{\mathbf U(\curlyu')}{\sigma}$ contains representation $\xi$ described below:
$$\mathbf U(\curlyu')/\mathbf U^1(\curlyu')\cong\GL_{n_1}(\kF)\times\ldots\times\GL_{n_k}(\kF)$$
where $n_1+\ldots+n_k=N$. Let $\xi_i$ be cuspidal representations of $\GL_{n_i}(\kF)$ and let $\xi$ be a lift of $\xi_1\otimes\ldots\otimes\xi_k$ to $\mathbf U(\curlyu')$. 

If $\xi_i\not\cong\xi_j$ for some $i$ and $j$, then $(\mathbf U(\curlyu'),\xi)$ is a split type of level $(0,0)$ and by \cite{bk}(8.4.1) we know that a supercuspidal representation cannot contain a split type. Hence $n_1=\ldots=n_k$ and $\xi_1\cong\ldots\cong\xi_k$. Then $(\mathbf U(\curlyu'),\xi)$ is a simple type, but it is not maximal, and by \cite{bk}(6.2.1) we know that supercuspidal representations contain only maximal simple types. 

So $\sigma$ is a lift of a cuspidal representation. Hence $(\uz,\sigma)$ is a maximal simple type occurring in $\pi$, so by \cite{bk}(6.2.3) $\pi\cong\cIndu{F^{\times}\uz}{G}{\tilde{\sigma}}$, where $\tilde{\sigma}$ is some extension of $\sigma$ to $F^{\times}\uz$. Hence, an irreducible representation $\sigma'$ of $\uz$ will occur in $\Rest{G}{\uz}{\pi}$ if and only if $\IG{\sigma'}{\sigma}{\uz}\neq\emptyset$. If $\sigma'$ is another irreducible summand of $\pi^{\uu{1}}$, then from above we get that $\sigma'$ is a lift of a cuspidal representation and $\sigma'$ intertwines with $\sigma$ in $G$. By \cite{bk}(5.7.1) there exists an $x\in G$, such that $\uz=\uz^x$ and $\sigma'\cong\sigma^x$, since $(\uz,\sigma')$ is a simple type. As $\curlyu$ is a maximal $\oF$ order in $A$, we get that $x\in F^{\times}\uz$, and hence $\sigma'\cong\sigma$. 

Therefore all irreducible factors of $\pi^{\uu{1}}$ will be isomorphic to $\sigma$, but Proposition \ref{exist} implies that $\sigma$ occurs in $\pi$ with multiplicity one. So $\pi^{\uu{1}}\cong\sigma$. 
\end{proof}
\begin{cor}\label{intzero} Let $\curlyu$ be a maximal hereditary $\oF$ order in $A$, $\sigma$ a lift of an irreducible cuspidal representation of $\uz/\uu{1}$ and let $\sigma'$ be a lift of any irreducible representation of $\uz/\uu{1}$, then $\sigma$ and $\sigma'$ intertwine in $G$ if and only if $\sigma'\cong\sigma$.
\end{cor}
\begin{proof} Apply Proposition \ref{levelzero} to $\cIndu{F^{\times}\uz}{G}{\tilde{\sigma}}$, where $\tilde{\sigma}$ is any extension of $\sigma$ to $F^{\times}\uz$.  
\end{proof}
\begin{cor}\label{intcusp} Let $\curlyu$ be a maximal hereditary $\oF$ order in $A$, $\sigma$ a lift of an irreducible cuspidal representation of $\uz/\uu{1}$ and let $\KK$ be a compact open subgroup of $\uz$, such that its image in $\uz/\uu{1}$ is a sufficiently small subgroup of $\uz/\uu{1}$. Moreover, if $q_F=2$ or $q_F=3$, then we assume further that the image of $\KK$ in $\uz/\uu{1}$ is contained in a proper parabolic subgroup of $\uz/\uu{1}$. Then for every irreducible summand $\xi$ of $\Rest{\uz}{\KK}{\sigma}$ there exists a lift $\sigma'$ of an irreducible representation of $\uz/\uu{1}$, such that $\brac{\xi}{\sigma'}{\KK}\neq 0$ and $\IG{\sigma}{\sigma'}{\uz}=\emptyset$.
\end{cor}  
\begin{proof} This is immediate from Proposition \ref{cuspidals} and Corollary \ref{intzero}.
\end{proof}
The following proposition can be easily obtained by making some cosmetic changes to \cite{bk}(5.3.2).
\begin{prop}\label{inter} Let $[\curlyu,n,0,\beta]$ be a simple stratum, $J=\JbU$, $J^1=\JibU{1}$ and $\theta\in\CUm{0}{\beta}$. Let $\eta$ be the unique representation of $J^1$ containing $\theta$ and let $\kappa$ be a $\beta$-extension of $\eta$. Let  $\zeta$ and $\zeta'$ be two lifts to $J$ of  irreducible representations of $J/J^1 \cong \ubz{\beta} /
\ub{1}{\beta}$. Suppose, that  $\IB{\beta}{\zeta}{\zeta'}{\ubz{\beta}}=\emptyset$, then $\IG{{\kappa\otimes\zeta}}{{\kappa\otimes\zeta'}}{J}=\emptyset$.
\end{prop}
We return to ideas and notations of Section \ref{important}.
\begin{prop}\label{doneA} Suppose $(g,Kg\kU)$ has property (A) and let  $\tau$ be an irreducible  representation of $K$, such that $\brac{\tau}{\Indu{{\tcap{g}{K}{\uz}}}{K}{\rho^g}}{K}\neq 0$, then $\tau$ cannot be a type.
\end{prop}
\begin{proof} Let $(J,\lambda)$ be a simple type, with the simple stratum $[\curlyu,n,0,\beta]$, such that $\rho\cong\Indu{J}{\uz}{\lambda}$ and
$\brac{\tau}{\Indu{{K\cap \cgate{g}{J}}}{K}{\lambda^g}}{K}\neq 0$. Let  $E=F[\beta]$. We have to consider two cases.

Suppose $\eUof=1$ and $(J,\lambda)=(\uz,\sigma)$, where $\sigma$ is a lift of a cuspidal representation of $\uz/\uu{1}$. Then $\Kg\le\uz$ and the image of $\Kg$ in $\uz/\uu{1}$ is a proper parabolic subgroup, so by Corollary \ref{intcusp} and Proposition \ref{key} $\tau$ cannot be a type.

Otherwise, $\lambda=\kappa\otimes\sigma$, where $\sigma$ is a lift of a cuspidal representation of $\ubz{\beta}/\ub{1}{\beta}$. Let $H$ be the image of $\Kg\cap J$ in $J/J^1\cong\ubz{\beta}/\ub{1}{\beta}$. By Corollary \ref{bananas} $H$ is a sufficiently small subgroup $\ubz{\beta}/\ub{1}{\beta}$. Moreover, if $q_E=2$ or $q_E=3$, then $E$ is a totally ramified extension of $F$, so $H$ is contained in a proper parabolic subgroup.

We will abuse the notation in the following way. Since $\uu{1}$ is a subgroup of $\Kg$, $$(J\cap\Kg)/J^1\cong(\ubz{\beta}\cap\Kg)/\ub{1}{\beta}$$ we will not distinguish between  representations of $J\cap\Kg$ (resp. $J$ ) on which $J^1$ acts trivially and representations of $\ubz{\beta}\cap\Kg$ (resp. $\ubz{\beta}$) on which $\ub{1}{\beta}$ acts trivially. 

Let $\xi$ be an irreducible summand of $\Rest{G}{J\cap\Kg}{\lambda}$, then $\brac{\xi}{\kappa\otimes\zeta}{J\cap\Kg}\neq0$, for some irreducible summand $\zeta$ of $\Rest{G}{\ubz{\beta}\cap\Kg}{\sigma}$. By Corollary \ref{intcusp} there exists a lift $\sigma'$ of an irreducible representation of $\ubz{\beta}/\ub{1}{\beta}$ to $\ubz{\beta}$ such that $\brac{\zeta}{\sigma'}{\ubz{\beta}\cap\Kg}\neq0$ and
$\sigma$ does not intertwine with $\sigma'$ in $B^{\times}_{\beta}$. Let $\lambda'=\kappa\otimes\sigma'$, then $\brac{\xi}{\lambda'}{J\cap\Kg}\neq0$ and 
by Proposition \ref{inter} $\lambda$ does not intertwine with $\lambda'$ in $G$. The representation $\lambda'$ is irreducible, since:
$$\dim\sigma'=\brac{\kappa\otimes\sigma'}{\eta}{J^1}=\brac{\kappa\otimes\sigma'}{\kappa\otimes\Indu{J^1}{J}{\Eins}}{J}=\sum_{\zeta'} \dim \zeta'\brac{\kappa\otimes\sigma'}{\kappa\otimes\zeta'}{J}$$
where $\eta=\Rest{G}{J^1}{\kappa}$ and the sum is taken over all the irreducible representations of $\ubz{\beta}/\ub{1}{\beta}$ lifted to $J$. The equality implies that $\brac{\lambda'}{\lambda'}{J}=1$. So by Proposition \ref{key} $\tau$ cannot be a type.
\end{proof}
\section{Double cosets with property (B)}\label{b}

If $(g,Kg\kU)$ has property (B) then the map $\uz\cap \cgate{{g^{-1}}}{K}\rightarrow\uz/\uu{1}$ is surjective, so we have to do something different than in the previous section. If $N=2$ the groups $\uu{1}\cap \cgate{{g^{-1}}}{K}$ coincide with the ones considered in \cite{ap}\S A.3.10. For general $N$ we will show that if $[\curlyu,n,0,\beta]$ is a simple stratum, then $\HibU{1}\cap\cgate{{g^{-1}}}{K}\neq\HibU{1}$ and will find $\theta'$ satisfying conditions of Proposition \ref{key}. Again it is more convenient to work with a larger subgroup than $\uu{1}\cap \cgate{{g^{-1}}}{K}$.

\begin{defi}\label{marvelous} Let $M$ be the following $\oF$-lattice in $A$: $$M=\{h\in \curlyp:hL_{e-1}\subseteq L_{e+1}\}$$where $e=\eUof$. And let $\KK=1+M$, be  a subgroup of $\uu{1}$.
\end{defi} 
Since $(g,Kg\kU)$ has property (B) $\uu{1}\cap \cgate{{g^{-1}}}{K}$ is a subgroup of $\KK$. Also, from the definition it is clear that $\uu{2}$ is a subgroup of $\KK$. 
\begin{lem}\label{KKK}Suppose $[\curlyu,n,0,\beta]$ is a simple stratum and $E=F[\beta]$, then $$\KK\cap\ub{1}{\beta}=1+M_{\beta}\mbox{, where } M_{\beta}=\{h\in\curlyq{\beta}: hL_{e_{\beta}-1}\subseteq L_{e_{\beta}+1}\}$$ and $e_{\beta}=\eBoe{\beta}$.
Moreover, if $\eBoe{\beta}=1$, then $\KK\cap\ub{1}{\beta}=\ub{2}{\beta}$.
\end{lem}
\begin{proof} If $x\in\KK\cap\ub{1}{\beta}$, then $x-1\in \curlyq{\beta}\cap M=\{h\in\curlyq{\beta}: hL_{e-1}\subseteq L_{e+1}\}$. Since $L_{i+me_{\beta}}=\pi_E^mL_{i}$, for all $i$, $e=e(E|F)e_{\beta}$ and  $x$ commutes with $\pi_E$, as $x\in B_{\beta}$, we have $x-1\in M_{\beta}$.

If $\eBoe{\beta}=1$, then $L_i=\pi^i_E L_0$, for all $i$, and since $x$ commutes with $\pi_E$, we have $(x-1)L_i\subseteq L_{i+2}$. That implies $x-1 \in \curlyqq{2}{\beta}$, so $x\in \ub{2}{\beta}$.
\end{proof}
We return to ideas and notations of Section \ref{important}.
\begin{prop}\label{nearlydoneB} Suppose $(g,Kg\kU)$ has property (B) and let  $\tau$ be an irreducible  representation of $K$, such that $\brac{\tau}{\Indu{{\tcap{g}{K}{\uz}}}{K}{\rho^g}}{K}\neq 0$. Moreover, let $(J,\lambda)$ be a simple type, with the simple stratum $[\curlyu,n,0,\beta]$,  such that $\rho\cong\Indu{J}{\uz}{\lambda}$ and $\brac{\tau}{\Indu{{K\cap \cgate{g}{J}}}{K}{\lambda^g}}{K}\neq 0$. Suppose $r=-k_0(\beta,\curlyu)>1$, then $\tau$ cannot be a type.
\end{prop}
\begin{proof} Let $E=F[\beta]$. Since $r>1$, \cite{bk}(3.1.15) implies the following decompositions:  $$\HibU{1}=\ub{1}{\beta}\HibU{2}$$    $$\HibU{1}\cap\KK=(\ub{1}{\beta}\cap\KK)\HibU{2}$$ as $\uu{2}$ is a subgroup of $\KK$. Now $(J,\lambda)$ is contained in a supercuspidal representation, so by \cite{bk}(6.2.1) $\eBoe{\beta}=1$. Hence by Lemma \ref{KKK} $$\KK\cap\HibU{1}=\HibU{2}$$
Let $\theta$ be a simple character occurring in $\Rest{J}{{H^1}}{\lambda}$. Since $\eBoe{\beta}=1$ the map $$H^1\rightarrow H^1/H^2\cong\ub{1}{\beta}/\ub{2}{\beta}\stackrel{\det_{B}}{\rightarrow}(1+\pE)/(1+\pE^2)$$ is surjective. Let $\tilde{\mu}$ be any non-trivial abelian character $$\tilde{\mu}:(1+\pE)/(1+\pE^2)\rightarrow\CC^{\times}$$
and $\mu$ be its lift to $H^1$. Let $\theta'=\theta\mu$, then $\Rest{H^1}{{H^1\cap\KK}}{\theta'}= \Rest{H^1}{{H^1\cap\KK}}{\theta}$.
 $$\IG{\theta'}{\theta}{H^1}\subseteq \IG{\theta'}{\theta}{H^2}=\IG{\theta}{\theta}{H^2}=J^1B_{\beta}^{\times}J^1$$
by \cite{bk}(3.3.2). By \cite{bk}(3.2.5) $\theta'$ is a simple character, so again by \cite{bk}(3.3.2)
 $$\IG{\theta'}{\theta'}{H^1}=J^1B_{\beta}^{\times}J^1\mbox{, and }\IG{\theta}{\theta}{H^1}=J^1B_{\beta}^{\times}J^1$$
that implies that $J^1B_{\beta}^{\times}J^1\subseteq \IG{\mu}{\mu}{H^1}$. Since all the representations above are $1$-dimensional, we may write everything explicitly. From above, if $\theta'$ and $\theta$ intertwine in $G$, then there exists $x\in J^1B_{\beta}^{\times}J^1$, such that $$\theta(h)\mu(h)=\theta(xhx^{-1}) \mbox{, }\every h\in H^1\cap x^{-1}H^1x$$ Since, such $x$ will intertwine $\theta$ with itself, we have $$\mu(h)=1\mbox{, } \every h\in H^1\cap x^{-1}H^1x$$ As $H^1$ is normal in $J^1$, we may assume $x=bj$, where $b\in B_{\beta}^{\times}$ and $j\in J^1$. Since $\mu(jhj^{-1})=\mu(h)$, for all $h\in H^1$, the intertwining of $\theta'$ and $\theta$ would imply that $$\mu(h)=1\mbox{, }\every h\in H^1\cap b^{-1}H^1b$$ By restricting to $\ub{1}{\beta}$, we get $$\ub{1}{\beta}\cap \cgate{b}{{\ub{1}{\beta}}}\le \Ker \mu$$ Since $\mu$ extends to $\ubz{\beta}$ and $\ub{1}{\beta}$ is normal in $\ubz{\beta}$, we have $$\ub{1}{\beta}\cap \cgate{b_1}{{\ub{1}{\beta}}}\le \Ker \mu\mbox{, }\every b_1\in \ubz{\beta}b\ubz{\beta}$$ We choose a basis, which identifies $\ubz{\beta}$ with $\GL_{\frac{N}{d}}(\oE)$, where $d=[E:F]$, and take $b_1$ to be a diagonal matrix with the eigenvalues equal to powers of $\pi_E$.  Conjugation by $b_1$ will fix the group $D$ of diagonal matrices in $\ub{1}{\beta}$ and $D\not\le\Ker \mu$. Hence $\theta$ and $\theta'$ do not intertwine in $G$. By Proposition \ref{key} $\tau$ cannot be a type.
\end{proof}
\begin{remar}When $N=2$ the arguments above are essentially \cite{ap}\S A.3.9 and \S A.3.10. 
\end{remar}
If $k_0(\beta,\curlyu)=-1$, then $\HibU{1}\neq\ub{1}{\beta}\HibU{2}$, and if $N=2$ Henniart uses a result of Casselman, which is not available for $N>2$, see \cite{ap}\S A.3.11. So we need a new idea. We recall some definitions. 

\begin{defi}\cite{bk}(2.3.1) A stratum of the form $[\curlyu,n,n-1,b]$ is called \textbf{fundamental} if $b+\curlypp{1-n}$ does not contain a nilpotent element
of $A$.
\end{defi}

Let $[\curlyu,n,n-1,b]$ be a fundamental stratum. We choose a prime element $\pi_F$ of $F$ and set $$y_b=b^{\frac{e}{m}}\pi_F^{\frac{n}{m}}+\curlyp,$$
 where $e=\eUof$, $m=\gcd(n,e)$. As an element of $\curlyu/\curlyp$, this depends only on the equivalence class of the stratum $[\curlyu,n,n-1,b]$.
 Let $\phi_b(X)\in\kF[X]$ be the characteristic polynomial of $y_b$ considered as an element of $\End_{\kF}(L_0/L_e)$ via the canonical embedding $\curlyu/\curlyp\subset\End_{\kF}(L_0/L_e)$.  
\begin{defi}\cite{bk}(2.3.3) A fundamental stratum  $[\curlyu,n,n-1,b]$ is called \textbf{split fundamental} if $\phi_b(X)$ has at least two distinct irreducible factors in $\kF[X]$. Otherwise, we say that $[\curlyu,n,n-1,b]$ is \textbf{non-split fundamental}.
\end{defi}
We start with the simplest case, when the simple stratum occurring in $\pi$ is $[\curlyu,1,0,\beta]$. The following Lemmas are preparation for Proposition \ref{onezero}. 

\begin{lem}\label{fund} Suppose $\curlyu$ is a principal hereditary $\oF$ order in $A$,
$e=\eUof$ and $b\in\curlypp{-1}$. We identify $\curlyu$ with block upper triangular matrices modulo $\pF$ and write  $b=(B_{ij})$, where $B_{ij}\in \MMe{\pF^{-1}}$, for $1\le i,j\le e$. Suppose $\pi_F B_{1e}, B_{21},\ldots,B_{(e-1)e}\in\GLNe$, then $[\curlyu,1,0,b]$ is a fundamental stratum.
Moreover, 
$$\phi_b(X)=(\det(X-\pi_F  B_{e(e-1)} \ldots B_{21} B_{1e}))^e \pmod{\pF}.$$ 
\end{lem} 
\begin{proof} Both statements above depend only on the coset $b+\curlyu$. We also know that $b\in\curlypp{-1}$ , so we may assume that $B_{1e}\in \MMe{\pF^{-1}}$, $B_{(i+1)i}\in \MMe{\oF}$, for $1\le i< e$, and $B_{ij}=0$, otherwise.

Let $\Pi$ be the element defined in the Section \ref{setup}. 
Then $\Pi b$ is a block diagonal matrix with the $i$-th block equal to $B_{(i+1)i}$ for $1\le i < e$ and the $e$-th  block equal to $\pi_F B_{1e}$. By our assumption, the blocks on diagonal are in $\GLNe$. So $\Pi b\in \uz$, hence 
$$b+\curlyu=\Pi^{-1}(\Pi b+\curlyp)\subset \Pi^{-1}\uz$$
as $\valu(\Pi)=1$. So $\valu(b)=-1$ and every element in $b+\curlyu$ is invertible, hence $[\curlyu,1,0,b]$ is a fundamental stratum.

Using block multiplication, we can calculate $\pi_F b^e$. Let $\pi_F b^e=(\tilde{B}_{ij})$, where $\tilde{B}_{ij}\in \MMe{\oF}$ for $1\le i,j\le e$, then $\tilde{B}_{ij}=0$, if $i\neq j$ and : 
$$\begin{aligned}
\tilde{B}_{11}&= \pi_F B_{1e} B_{e(e-1)} B_{(e-1)(e-2)} \ldots B_{21},\\
\tilde{B}_{22}&=\pi_F B_{21}B_{1e} B_{e(e-1)} B_{(e-1)(e-2)} \ldots B_{32},\\
\vdots \\
\tilde{B}_{ee}&=\pi_F  B_{e(e-1)} B_{(e-1)(e-2)} \ldots B_{21} B_{1e}.
\end{aligned}$$ 
Hence $$\phi_b(X)=(\det(X-\pi_F  B_{e(e-1)} \ldots B_{21} B_{1e}))^e \pmod{\pF}.$$
\end{proof}

\begin{lem}\label{eNtwo} Let $\curlyu$ be a principal hereditary $\oF$ order in $A$, $q_F=2$ and $\eUof=\frac{N}{2}$. Suppose the stratum $[\curlyu,1,0,b]$ is fundamental and $\phi_b(X)$ is a power of $X^2+X+1$, then
$[\curlyu,1,0,b]$ is equivalent to a simple stratum.
\end{lem}
\begin{proof} The polynomial $X^2+X+1$ is irreducible over $\FF_2$, so the stratum $[\curlyu,1,0,b]$ is non-split fundamental. From \cite{bk}(2.3.4), we know that there exists a simple stratum $[\curlyu',n',n'-1, \alpha]$ such that $b+\curlyu\subseteq\alpha+\curlyp{'}^{1-n'}$, moreover $\frac{2}{N}=\frac{n'}{e(\curlyu')}$, and the lattice chain defining $\curlyu'$ contains that defining $\curlyu$.

If $e(\curlyu'|\oF)=\frac{N}{2}$, then $\curlyu=\curlyu'$ and $n'=1$. So $b+\curlyu=\alpha+\curlyu$ and we are done.

Otherwise, $e(\curlyu'|\oF)=N$, $n'=2$ and $b+\curlyp{'}^{-1}=\alpha+\curlyp{'}^{-1}$. Hence $\nu_{\curlyu'}(b)=-2$, so $\pi_F b^{\frac{N}{2}}\in \curlyu'$. Therefore the characteristic polynomial of $\pi_F b^{\frac{N}{2}}$ modulo $\pF$ is a power of $X-1$, but  $\phi_b(X)$ associated to $[\curlyu,1,0,b]$ is also the characteristic polynomial of $\pi_F b^{\frac{N}{2}}$ modulo $\pF$ and it is a power of $X^2+X+1$. We get a contradiction. 
\end{proof}

\begin{lem}\label{onezerosimple} Suppose $[\curlyu,1,0,\beta]$ and $[\curlyu,1,0,\gamma]$ are simple strata, such that $\IG{\psb{\beta}}{\psb{\gamma}}{\uu{1}}\neq \emptyset$, then $\phi_{\beta}(X)=\phi_{\gamma}(X)$. 
\end{lem}
\begin{proof} In this case we have $\HitU{1}=\HigU{1}=\uu{1}$, $\psb{\beta}\in\CUm{0}{\beta}$ and $\psb{\gamma}\in\CUm{0}{\gamma}$. As $\psb{\beta}$ and $\psb{\gamma}$ intertwine in $G$ we use \cite{bk}(3.5.11) to  get $x\in \uz$ such that $\psb{\beta}=\psb{\gamma}^x$. That implies $\beta+\curlyu=x^{-1}\gamma x+\curlyu$, and as conjugation does not change  characteristic polynomials we get the result.  
\end{proof}
\begin{defi} For an $\oF$ lattice $L$ in $A$, we define $$L^*=\{x\in A: \psa(xh)=1\mbox{ for all }h \in L\}$$
\end{defi}

\begin{prop}\label{onezero} Let $[\curlyu,1,0,\beta]$ be a simple stratum, $e=\eUof>1$, $E=F[\beta]$, $\eBoe{\beta}=1$ and $M$ as in Definition \ref{marvelous}. Then there exists $b\in M^*$ such that one of the following holds:
\begin{enumerate}
\item If $q_F>2$ and $e<N$ or $q_F=2$ and $e<\frac{N}{2}$ then the stratum $[\curlyu,1,0,\beta+b]$ is split fundamental.
\item If $q_F=2$, $e=\frac{N}{2}$ and $E$ is totally ramified over $F$ then $[\curlyu,1,0,\beta+b]$ is equivalent to a simple stratum and $\IG{\psb{\beta}}{\psb{\beta+b}}{\uu{1}}=\emptyset$.
\item If $e=N$ or $q_F=2$, $e=\frac{N}{2}$ and $E$ is not totally ramified over $F$ then $[\curlyu,1,0,\beta+b]$ is not fundamental and  $\IG{\psb{\beta}}{\psb{\beta+b}}{\uu{1}}=\emptyset$.
\end{enumerate}
\end{prop}
\begin{remar}Note, that $b\in M^*$ if and only if $\Rest{G}{\KK}{\psb{\beta}}=\Rest{G}{\KK}{\psb{\beta+b}}$.
\end{remar}
\begin{proof} We identify $\curlyu$ with block upper triangular matrices modulo $\pF$. Since $\beta\in\curlyp^{-1}$, we can write $\beta$ with respect to our fixed basis as a matrix $(A_{ij})$, $1\le i,j \le e$, where $A_{1e}\in \MMe{\pF^{-1}}$, $A_{ij}\in \MMe{\oF}$ for all $i\le j+1$ and $(i,j)\neq(1,e)$ and $A_{ij}\in \MMe{\pF}$ otherwise. Let $y$ be the matrix $(\tilde{A}_{ij})$, where $\tilde{A}_{(i+1)i}=A_{(i+1)i}$, for $1\le i<e$, $\tilde{A}_{1e}=A_{1e}$, and $\tilde{A}_{ij}=0$, otherwise. So $y+\curlyu=\beta+\curlyu$. Let $\Pi$ be the element defined in the Section \ref{setup}. Since $\valu(\Pi)=1$ we have 
$$\Pi y\in \Pi(\beta+\curlyu)=\Pi\beta{\uu{1}}\subset\uz.$$
The matrix $\Pi y$ is block diagonal with the $i$-th block equal to $A_{(i+1)i}$ for $1\le i < e$ and the $e$-th  block equal to $\pi_F A_{1e}$. So we can apply
Lemma \ref{fund} to get that stratum $[\curlyu,1,0,\beta]$ is fundamental and 
$$\phi_{\beta}(X)=\phi_y(X)=(\det(X-\pi_F  A_{e(e-1)} \ldots A_{21} A_{1e}))^e \pmod{\pF}.$$ 
Let $b=(B_{ij})$, $1\le i,j\le e$, where $B_{1e}\in \MMe{\pF^{-1}}$ and $B_{ij}=0$ otherwise. Then from multiplication of blocks, we can see that $b\in M^*$. In each case we will find a matrix $B_{1e}$, such that conditions of proposition are satisfied. We note that, $\eBoe{\beta}=1$ implies that the ramification index of $E/F$ $e(E|F)=\eUof$.

1. If $q_F>2$ and $e<N$ or $q_F=2$ and $e<\frac{N}{2}$, then we can find a matrix $C\in \GLNe$, such that the characteristic polynomial of $C$ in modulo $\pF$ contains two distinct irreducible factors over $\kF$. Note that, this is not possible if $q_F=2$ and $e=\frac{N}{2}$. Let
$$B_{1e}= \pi_F^{-1}  A_{21}^{-1} \ldots A_{e(e-1)}^{-1}C-A_{1e}$$
then Lemma \ref{fund} implies that the stratum $[\curlyu,1,0,\beta+b]$ is fundamental and $\phi_{\beta+b}(X)=(\det(X-C))^e\pmod{\pF}$, so $[\curlyu,1,0,\beta+b]$ is split fundamental.

2. If $q_F=2$, $e=\frac{N}{2}$ and $E$ is totally ramified over $F$, let $C\in\MM(2,\oF)$ be a matrix such that the characteristic polynomial of $C$  modulo $\pF$ is $X^2+X+1$.
Let $$B_{1e}= \pi_F^{-1}  A_{21}^{-1} \ldots A_{e(e-1)}^{-1}C-A_{1e}$$
so $[\curlyu,1,0,\beta+b]$ is fundamental, $\phi_{\beta+b}(X)=(\det(X-C))^e\pmod{\pF}$, which is a power of $X^2+X+1$, so by Lemma \ref{eNtwo} $[\curlyu,1,0,\beta+b]$ is equivalent to a simple stratum.
As $[\curlyu,1,0,\beta]$ is simple, we have $\valu(\beta)=k_0(\beta,\curlyu)=-1$ and by \cite{bk}(1.4.15) $\pi_F\beta^{\frac{N}{2}}+\pE$ generates $\kE$ over $\kF$. As $E$ is totally ramified over $F$, we get that $\phi_{\beta}(X)$ is a power of $X-1$. Then Lemma \ref{onezerosimple} implies that $\psi_{\beta}$ and $\psi_{\beta+b}$ do not intertwine in $G$.

3. If $e=N$ or $q_F=2$, $e=\frac{N}{2}$ and $E$ is not totally ramified over $F$, then $[E:F]=N$. Let $J=J(\beta,\curlyu)$, then $J^1(\beta)=H^1(\beta)=\uu{1}$, $J/\uu{1}\cong\kE^{\times}$ and $\psi_{\beta}\in \CUm{0}{\beta}$. Let $\lambda$ be a simple type, such that $$\Rest{J}{\uu{1}}{\lambda}\cong\psi_{\beta}$$ and $\Lambda$ any extension of $\lambda$ to $E^{\times}J$. Let $$\pi'=\cIndu{E^{\times}J}{G}{\Lambda}$$ so $\pi'$ is a supercuspidal representation. As $\uu{1}$ is the unique maximal pro-$p$ subgroup of $J$, we have:
$$\Rest{G}{\uu{1}}{\pi'}\cong \Mplus{g}{\uu{1}}{G}{E^{\times}J}\IndResg{g}{\uu{1}}{G}{\uu{1}}{\psi_{\beta}}$$  Let $B_{1e}=-A_{1e}$, then the stratum $[\curlyu,1,0,\beta+b]$ is not fundamental. Suppose $\psi_{\beta+b}$ and $\psi_{\beta}$ intertwine in $G$, then from above we know that $\psi_{\beta+b}$ occurs in $\Rest{G}{\uu{1}}{\pi'}$.

Let $\curlyu_M=\End_{\oF}(L_0)$, then $\curlyu_M$ is a maximal hereditary $\oF$ order in $A$, $K=\mathbf U(\curlyu_M)$ and since $\mathbf U^1(\curlyu_M)=I_N+\MM(N,\pF)$, where $I_N$ is the identity matrix, we have $\Rest{\uu{1}}{ \mathbf U^1(\curlyu_M)}{\psi_{\beta+b}}=1$. 

So $\pi{'}^{\mathbf U^1(\curlyu_M)}\neq 0$, hence Proposition \ref{levelzero} implies that $\Rest{G}{ \mathbf U(\curlyu_M)}{\pi'}$ contains $\sigma$, which is a lift of a cuspidal representation of  $ \mathbf U(\curlyu_M)/\mathbf U^1(\curlyu_M)$.

Since $(\mathbf U(\curlyu_M),\sigma)$ is another simple type occurring in $\pi'$, by \cite{bk}(6.2.4) there exists $g\in G$, such that $\mathbf U(\curlyu_M)=J^g$ and $\sigma\cong\lambda^g$. But $J$ has a unique maximal pro-$p$ subgroup, and $\mathbf U(\curlyu_M)$ does not, so that cannot happen. We get a contradiction, so $\psi_{\beta+b}$ does not intertwine with $\psi_{\beta}$.  
\end{proof} 

We recall the following definition.
\begin{defi}\label{splittwo}\cite{bk}(8.1.1) A \textbf{split type of level $(x,x)$, $x>0$}, is a pair $(K',\vartheta)$ given as follows:
\begin{itemize}
 \item[(i)] $[\curlyu,n,n-1,b]$ is a split fundamental stratum in $A$
 \item[(ii)] $n>0$, $\gcd(n,e(\curlyu))=1$, $x=n/e(\curlyu)$
 \item[(iii)] $K'=\uu{n}$, $\vartheta=\psb{b}$.
\end{itemize}
\end{defi}  

\begin{cor}\label{almostdoneB} Suppose $(g,Kg\kU)$ has property (B) and let  $\tau$ be an irreducible  representation of $K$, such that $\brac{\tau}{\Indu{{\tcap{g}{K}{\uz}}}{K}{\rho^g}}{K}\neq 0$.

Moreover, suppose $(J,\lambda)$ is a simple type with the simple stratum $[\curlyu,1,0,\beta]$, such that $\rho\cong\Indu{J}{\uz}{\lambda}$ and $\brac{\tau}{\Indu{{K\cap \cgate{g}{J}}}{K}{\lambda^g}}{K}\neq 0$, then $\tau$ cannot be a type.
\end{cor}  
\begin{proof} Since $[\curlyu,1,0,\beta]$ is a simple stratum, we have $r=-k_0(\beta,\curlyu)=1$ and since $(J,\lambda)$ is a simple type occurring in a supercuspidal representation, we have $\eBoe{\beta}=1$. Also from the definitions of $\JibU{1}$ and $\HibU{1}$ \cite{bk}(3.1.7) and (3.1.8) we get $$\JibU{1}=\HibU{1}=\uu{1}$$ and from the definition of simple characters \cite{bk}(3.2.1) we have that the only simple character is $\psb{\beta}$. Apply Proposition \ref{onezero} to this situation, then we get  an abelian character $\psb{\beta+b}$ of $\uu{1}$, such that $$\Rest{\uu{1}}{\KK}{\psb{\beta+b}}=\Rest{\uu{1}}{\KK}{\psb{\beta}}.$$ From the proof of Proposition \ref{key}, one can see that if $\tau$ is a type, then $\Rest{G}{\uu{1}}{\pi}$ must contain all irreducible summands of $\Indu{\uu{1}\cap K^{g^{-1}}}{\uu{1}}{\tau^{g^{-1}}}$, so $\psb{\beta+b}$ must occur in $\Rest{G}{\uu{1}}{\pi}$. 

If $q_F>2$ and $\eUof<N$ or $q_F=2$ and $\eUof<\frac{N}{2}$, then $[\curlyu,1,0,\beta+b]$ is a split fundamental stratum, so  $(\uu{1},\psb{\beta+b})$ is a split type of level $(1/e,1/e)$. But by \cite{bk}(8.4.1) a supercuspidal representation cannot contain a split type. So $\tau$ is not a type.

In all the other cases of Proposition \ref{onezero} $\IG{\psb{\beta}}{\psb{\beta+b}}{\uu{1}}=\emptyset$. So we apply Proposition \ref{key} with $\theta=\psb{\beta}$ and $\theta'=\psb{\beta+b}$, and hence $\tau$ cannot be a type.
\end{proof}
We recall the following definition.
\begin{defi}\label{splitthree}\cite{bk}(8.1.3) A \textbf{split type of level $(x,y)$, $x>y>0$}, is a pair $(K',\vartheta)$ given as follows:
\begin{itemize}
\item[(i)] $[\curlyu,n,m,\beta]$ is a simple stratum in $A$ with $E=F[\beta]$, $B=\End_E(V)$, $\urlyb=\curlyu\cap B$, $e_{\beta}=e(\urlyb|\oE)$, $\gcd(m,e_{\beta})=1$, $x=n/e(\curlyu)$, $y=m/e(\curlyu)$
\item[(ii)] $K'=H^m(\beta,\curlyu)$
\item[(iii)] $\vartheta=\theta \psb{c}$, for some $\theta\in\CUm{m-1}{\beta}$ and some $c\in \curlyp^{-m}$, such that the stratum  $[\urlyb,m,m-1,s_{\beta}(c)]$ is split fundamental, where $s_{\beta}$ denotes a tame corestriction on $A$ relative to $E/F$.
\end{itemize}
\end{defi}   
 
\begin{prop}\label{doneB} Suppose $(g,Kg\kU)$ has property (B) and let  $\tau$ be an irreducible  representation of $K$, such that $\brac{\tau}{\Indu{{\tcap{g}{K}{\uz}}}{K}{\rho^g}}{K}\neq 0$.

Moreover, let $(J,\lambda)$ be a simple type with the simple stratum $[\curlyu,n,0,\beta]$, such that  $\rho\cong\Indu{J}{\uz}{\lambda}$, $\brac{\tau}{\Indu{{K\cap \cgate{g}{J}}}{K}{\lambda^g}}{K}\neq 0$. Suppose that  $r=-k_0(\beta,\curlyu)=1$ and $n>1$, then $\tau$ cannot be a type.
\end{prop}  
\begin{proof} Let $\theta\in\CUm{0}{\beta}$, such that $\theta$ occurs in $\Rest{\JbU}{\HitU{1}}{\lambda}$. By \cite{bk}(3.2.3) there exists a simple stratum $[\curlyu,n,1,\gamma]$, such that $[\curlyu,n,1,\beta]\sim [\curlyu,n,1,\gamma]$, $\HibU{1}=H^1(\gamma,\curlyu)$ and   $$\theta=\theta_0 \psb{c},$$ where $\theta_0\in\CUm{0}{\gamma}$ and $c=\beta-\gamma$. Since $\beta+\curlyp^{-1}=\gamma+\curlyp^{-1}$, we have $\valu(c)\ge -1$. If $\valu(c)\ge 0$, we would have $\beta+\curlyu=\gamma+\curlyu$, so $[\curlyu,n,0,\beta]\sim [\curlyu,n,0,\gamma]$. Since $[\curlyu,n,0,\beta]$ and $[\curlyu,n,0,\gamma]$ are both simple \cite{bk}(2.4.1)(ii)(a) would imply $k_0(\beta,\curlyu)=k_0(\gamma,\curlyu)$, but since $[\curlyu,n,1,\gamma]$ is a simple stratum, we have $k_0(\gamma,\curlyu)\le -2$ and $k_0(\beta,\curlyu)=-1$ , hence $$\valu(c)=-1.$$
That implies $\psb{c}$ extends to an abelian character of $\uu{1}$ and $\Rest{G}{\uu{2}}{\psb{c}}=1$. 

Since $k_0(\gamma,\curlyu)\le -2$, we have $$\HigU{1}=\ub{1}{\gamma}\HigU{2}$$ and $\uu{2}$ is a subgroup of  $\KK$, so $$\KK\cap\HitU{1}=\KK\cap\HigU{1}=(\ub{1}{\gamma}\cap\KK)\HigU{2}=\KK_{\gamma}\HigU{2}$$
where $\KK_{\gamma}=1+M_{\gamma}=\ub{1}{\gamma}\cap\KK$ as in Lemma \ref{KKK}.  
Let $$e_{\gamma}=e(\urlyb_{\gamma}|\oo_{F[\gamma]})$$
If $e_{\gamma}=1$, then by Lemma \ref{KKK} $\KK\cap \HigU{1}=\HigU{2}$. Let $\theta'=\theta_0$, then $\Rest{H^1}{\HigU{1}\cap \KK}{\theta'}=\Rest{H^1}{\HigU{1}\cap \KK}{\theta}$ and by \cite{bk}(3.5.12) $\IG{\theta}{\theta'}{\HibU{1}}=\emptyset$. So by Proposition \ref{key} $\tau$ cannot be a type.

If $e_{\gamma}>1$, we fix a continuous character $\psi_{F[\gamma]}$ of the additive group $F[\gamma]$ with the conductor $\pp_{F[\gamma]}$ and let $$\psi_{B_{\gamma}}(b')=\psi_{F[\gamma]}(\tr_{B_{\gamma}/F[\gamma]}(b'))\mbox{, } \every b'\in B_{\gamma}$$ Then there exists a tame corestriction $s_{\gamma}$   on $A$ relative to $F[\gamma]/F$, such that $$\psi_{A}(ab')=\psi_{B_{\gamma}}(s_{\gamma}(a)b') \mbox{, }\every a\in A \mbox{, }\every b'\in B_{\gamma}$$
In particular, for every $c'\in \curlyp^{-1}$ we have
$$\psb{c',A}(b')=\psb{s_{\gamma}(c'),B_{\gamma}}(b')\mbox{, }\every b'\in \ub{1}{\gamma}$$
By \cite{bk}(2.4.1)(iii) there exists a simple stratum $[\urlyb_{\gamma},1,0, \delta]$ in $B_{\gamma}$, such that $$[\urlyb_{\gamma},1,0, s_{\gamma}(c)]\sim [\urlyb_{\gamma},1,0, \delta]$$
We want to apply the Proposition \ref{onezero} to $[\urlyb_{\gamma},1,0, \delta]$. Let $B_{\gamma,\delta}$ be the $B_{\gamma}$-centraliser of $F[\gamma,\delta]$ and $\curlyb{\gamma,\delta}=\curlyb{\gamma}\cap B_{\gamma,\delta}$. 

We claim that $e(\curlyb{\gamma,\delta}|\oo_{F[\gamma,\delta]})=1$. By \cite{bk}(2.2.8) we have $$e(F[\gamma,\delta]|F)=e(F[\beta]|F).$$
Since $e(\curlyb{\beta}|\oo_{F[\beta]})=1$, we also have $$\eUof=e(F[\beta]|F).$$
And $$\eUof=e(\curlyb{\gamma}|\oo_{F[\gamma]})e(F[\gamma]|F).$$
Hence $$e(F[\gamma,\delta]|F[\gamma])=e(\curlyb{\gamma}|\oo_{F[\gamma]}),$$
which proves the claim. So  we can apply Proposition \ref{onezero} to $[\urlyb_{\gamma},1,0, \delta]$. We get $d\in\curlyqq{-1}{\gamma}$, such that $\psb{\delta+d}$ is an abelian character of $\ub{1}{\gamma}$ and  
$$\Rest{\ub{1}{\gamma}}{\KK_{\gamma}}{\psb{\delta+d}}=\Rest{\ub{1}{\gamma}}
{\KK_{\gamma}}{\psb{\delta}}.$$
 By \cite{bk}(1.4.7) $s_{\gamma}:\curlypp{-1}\rightarrow\curlyqq{-1}{\gamma}$ is surjective. Choose $b\in\curlypp{-1}$, such that $s_{\gamma}(b)=d$, and let $\theta'=\theta_0\psb{c+b}$. If $h\in \KK\cap\HigU{1}$, then $h= h_1 h_2$, for some $h_1\in\KK_{\gamma}$, $h_2\in \HigU{2}$, and
$$\psb{c+b,A}(h)=\psb{c+b,A}(h_1)=\psb{s_{\gamma}(c+b),B_{\gamma}}(h_1)=\psb{\delta+d,B_{\gamma}}(h_1)$$
$$\psb{c,A}(h)=\psb{c,A}(h_1)=\psb{s_{\gamma}(c),B_{\gamma}}(h_1)=\psb{\delta,B_{\gamma}}(h_1)$$
From above  $\Rest{H}{\KK\cap\HigU{1}}{\psb{c+b}}=\Rest{H}{\KK\cap\HigU{1}}{\psb{c}}$ and hence $\Rest{H}{\KK\cap\HigU{1}}{\theta'}=\Rest{H}{\KK\cap\HigU{1}}{\theta}$. So if $\tau$ was a type, then by arguments in Proposition \ref{key}, we would have that $\theta'$ occurs in $\Rest{G}{\HitU{1}}{\pi}$. 

Suppose $q_{F[\gamma]}>2$ and $e_{\gamma}[F[\gamma]:F]<N$ or $q_{F[\gamma]}=2$ and $2e_{\gamma}[F[\gamma]:F]<N$, then the stratum $[\urlyb_{\gamma},1,0,s_{\gamma}(c+b)]$ is split fundamental, so $(H^1(\gamma,\curlyu), \theta')$ is a split type of level $(n/e,1/e)$, and by \cite{bk}(8.4.1), a supercuspidal representation cannot contain a split type. So $\tau$ cannot be a type.
 
In all the other cases of Proposition \ref{onezero}, $\psb{\delta}$ and $\psb{\delta+d}$ do not intertwine in $B^{\times}_{\gamma}$. We will show  that this implies that $\theta$ and $\theta'$ do not intertwine in $G$. 
$$\IG{\theta'}{\theta}{\HigU{1}}\subseteq\IG{\theta'}{\theta}{\HigU{2}}=\IG{\theta_0}{\theta_0}{\HigU{2}}=J^1(\gamma)B^{\times}_{\gamma}J^1(\gamma)$$
by \cite{bk}(3.3.2). By the same theorem $\IG{\theta_0}{\theta_0}{\HigU{1}}=J^1(\gamma)B^{\times}_{\gamma}J^1(\gamma)$ and $\theta_0$ is an abelian character, so if $h$ intertwines $\theta$ and $\theta'$ in $G$, it must also intertwine $\psb{c}$ and $\psb{c+b}$. Both characters extend to $\uu{1}$ and $\HigU{1}$ is normal in $J^1({\gamma})$, so if $h=j_1b'j_2$, where $j_1,j_2\in J^1({\gamma})$ and $b'\in B^{\times}_{\gamma}$, then $b'$ must also intertwine $\psb{c}$ and $\psb{c+b}$ in $G$ and hence $b'$ must intertwine the restrictions of these characters to $\ub{1}{\gamma}$ in $B^{\times}_{\gamma}$. So $$b'\in\IB{\gamma}{\psb{c}}{\psb{c+b}}{\ub{1}{\gamma}}=\IB{\gamma}{\psb{\delta}}{\psb{\delta+d}}{\ub{1}{\gamma}}=\emptyset$$
That implies $\IG{\theta'}{\theta}{\HibU{1}}=\emptyset$. By Proposition \ref{key} $\tau$ cannot be a type.
\end{proof}
\begin{remar}If $N=2$, the case above does not have to be considered. Since, we can always find a smooth quasicharacter $\chi$ of $F^{\times}$, such that the simple stratum $[\curlyu,n,0,\beta]$ occurring in $\pi\otimes\chi\circ\det$ has $\beta$ minimal over $F$,i.e., $\valu(\beta)=k_0(\beta,\curlyu)$. Then it is easy to see, that it is enough to prove the unicity of types for $\pi\otimes\chi\circ\det$. I was told by Bushnell, that this works if and only if $N$ is prime.
\end{remar}
\section{Inertial correspondence}\label{over}
We collect all the bits together.
\begin{thm}(Main)\label{nearly} Let $G=\GL_N(F)$ and $\pi$ be a smooth irreducible supercuspidal representation of $G$, then there exists a unique (up to isomorphism) smooth irreducible representation $\tau$ of $K=\GL_N(\oF)$, such that for any infinite dimensional smooth irreducible representation $\pi'$ of $G$:
 $$\Rest{G}{K}{\pi'}\mbox{ contains }\tau \Leftrightarrow \pi'\cong \pi \otimes \chi \circ \det $$
where $\chi$ is some unramified quasicharacter of $F^{\times}$.

Moreover, if $(J,\lambda)$ is a simple type in a sense of \cite{bk}, with the simple stratum $[\curlyu,n,0,\beta]$, such that $\uz\le K$ and  $\pi\cong\cIndu{{E^{\times}J}}{G}{\Lambda}$, where $E=F[\beta]$ and $\Lambda$ is an extension of $\lambda$ to $E^{\times}J$, then $\tau\cong\Indu{J}{K}{\lambda}$.

Further, $\tau$ occurs in $\Rest{G}{K}{\pi}$ with multiplicity one. 
\end{thm}
\begin{proof}Let $\tau$ be any irreducible representation of $K$ occurring in $\Rest{G}{K}{\pi}$.  
$$\Rest{G}{K}{\pi}\cong \Mplus{g}{K}{G}{{\kU}}\IndResg{g}{K}{G}{\uz}{\rho}$$
where $\rho=\Indu{J}{\uz}{\lambda}$. 
Hence $\brac{\tau}{\Indu{{\tcap{g}{K}{\uz}}}{K}{\rho^g}}{K}\neq 0$, for some representative $g\in G$. 

If the double coset $Kg\kU=K\kU$, then Proposition \ref{exist} says that $\tau\cong\Indu{J}{K}{\lambda}$, is a type and occurs in $\Rest{G}{K}{\pi}$ with multiplicity one.  

If the double coset $Kg\kU\neq K\kU$, then we combine Propositions \ref{doneA}, \ref{nearlydoneB}, \ref{almostdoneB} and \ref{doneB}, to get that $\tau$ cannot be a type. That establishes uniqueness. 
\end{proof}

It also allows us to define a kind of inertial local Langlands correspondence for supercuspidals.
\begin{cor}Let $W_F$ be the Weil group of $F$, $I_F$ the inertia subgroup,  $\varphi$ be a smooth $N$-dimensional representation of $I_F$, such that it extends to a smooth \underline{irreducible} Frobenius semisimple representation of $W_F$, then there exists a unique (up to isomorphism) smooth irreducible representation $\tau(\varphi)$ of $K=\GL_N(\oF)$, such that for any smooth irreducible infinite dimensional representation $\pi'$ of $G=\GL_N(F)$, $\tau(\varphi)$ occurs in $\pi'$ with multiplicity at most $1$ and :
$$\Rest{G}{K}{\pi'}\mbox{ contains }\tau(\varphi)\Leftrightarrow \Rest{{W_F}}{{I_F}}{\WD(\pi)}\cong\varphi$$
where $\WD(\pi)$ is a Weil-Deligne representation of $W_F$ corresponding to $\pi'$ via the local Langlands correspondence.
\end{cor}
\begin{proof} Let $\varphi_1$ be an irreducible smooth Frobenius semisimple representation of $W_F$, such that $\Rest{W_F}{I_F}{\varphi_1}\cong\varphi$ and  let $\mathcal L$ denote the Langlands correspondence going from the Galois side to the automorphic side. Local Langlands correspondence preserves tensoring with quasicharacters, and irreducible $N$-dimensional representations of $W_F$ are mapped to supercuspidal representations of $G$. So $\LL{\varphi_1}$ is supercuspidal and if $\pi'\in \mathfrak I(\LL{\varphi_1})$, then $\Rest{G}{I_F}{WD(\pi')}\cong\varphi$. Conversely, if $\varphi_2\cong\varphi_1\otimes \chi$, then $\LL{\varphi_2}\in \mathfrak I(\LL{\varphi_1})$. So it is enough to prove the following statement:

Let $\varphi_2$ be a smooth Frobenius semisimple representations of $W_F$, such that $\Rest{W_F}{I_F}{\varphi_2}\cong\varphi$, then
$\varphi_2\cong\varphi_1\otimes \chi $, where $\chi$ is some unramified quasicharacter of $F^{\times}$. 

Then Theorem \ref{nearly} applied to $\mathfrak I(\LL{\varphi_1})$ provides us with the unique $\tau=\tau(\varphi)$.

By tensoring with some unramified quasicharacter, we may assume that the image of $\varphi_1(W_F)$ in $\GL_N(\CC)$ is finite. By tensoring each irreducible factor of $\varphi_2$ by an unramified quasicharacter, we may assume that the image of $\varphi_2(W_F)$ in $\GL_N(\CC)$ is  also finite. We can view $\varphi_1$ and $\varphi_2$ as representations of a finite group $H=W_F/(\Ker\varphi_1\cap\Ker\varphi_2)$, and let $I$ be the image of inertia in $H$. Then $$0\neq\brac{\varphi_2}{\varphi}{I}=\brac{\varphi_2}{\Indu{I}{H}{\varphi}}{H}=
 \brac{\varphi_2}{\varphi_1\otimes\Indu{I}{H}{\Eins}}{H}$$ 
Since $I$ is normal in $H$ and $H/I$ is cyclic, we have $\brac{\varphi_2}{\varphi_1\otimes\chi}{H}\neq 0$, for some  $\chi$  an abelian character of $H/I$. Since $\varphi_1$ is irreducible and has the same dimension as $\varphi_2$, we get $\varphi_2\cong\varphi_1\otimes\chi$. \end{proof}

\end{document}